\documentclass[a4paper,12pt]{article}

\usepackage{amsmath}
\usepackage{a4wide}
\usepackage{amsfonts}
\newenvironment{pf}{\textbf{Proof:}}{\hspace{\stretch{1}}$\square$}

\newtheorem{thm}{Theorem}
\newtheorem{df}{Definition}
\newtheorem{pr}{Proposition}
\newtheorem{lem}{Lemma}
\newtheorem{dthm}{Definition-Theorem}

\title{Existence, Uniqueness and Approximation of a Stochastic Schr\"odinger Equation: the Poisson Case}
\author{Cl\'ement PELLEGRINI\\
\scriptsize{Institut C.Jordan}\\
\scriptsize{Universit\'e C.Bernard, lyon 1}\\
\scriptsize{21, av Claude Bernard}\\
\scriptsize{69622 Villeurbanne Cedex}\\
\scriptsize{France}\\
\scriptsize{e-mail: pelleg@math.univ-lyon1.fr}}
\begin{document}
\maketitle

\begin{abstract}
In quantum physics, recent investigations deal with the so-called "quantum trajectory" theory. Heuristic rules are usually used to
give rise to ``stochastic Schr\"odinger equations'' which are stochastic differential equations of non-usual type describing the
physical models. These equations pose tedious problems in terms of mathematical and physical justifications: notion of solution,
existence, uniqueness...

In this article, we concentrate on a particular case: the Poisson case. Random Measure theory is used in order to give rigorous sense
to such equations. We prove existence and uniqueness of a solution for the associated stochastic equation. Furthermore, the stochastic
model is physically justified by proving that the solution can be obtained as a limit of a concrete discrete time physical model.

\end{abstract}
\section*{Introduction}

Many recent developments in quantum mechanics deal with \textit{``Stochastic
Schr\"odinger Equations''}. These equations are classical stochastic
differential equations (also called \textit{"Belavkin equations"}) which
describe random phenomena in continuous measurement theory. The
solutions of these equations are called \textit{``quantum trajectories''}, they give account of the time evolution of an open quantum system undergoing
continuous measurement.

 Usually, in Quantum Optics or Quantum Communication, indirect
 measurement is performed in order to avoid phenomena like \textit{Zenon effect}. The physical
 setup is the one of a small system (an open system) interacting with an environment, the measurement is then performed on the environment.
 In the literature, stochastic Schr\"odinger equations are expressed as perturbations of the Master Equation
 which describes normally the evolution of the small system without measurement.

 Belavkin equations are stochastic Schr\"odinger equations describing the evolution of qubit, i.e two level system. Essentially there exists two type of equations.
\begin{itemize}
\item A diffusive equation
\begin{equation}\label{difffff}
d\rho_t=L(\rho_t)dt+\big[\rho_t C^\star+ C\rho_t-Tr\left[\rho_t(C+C^\star) \right]\rho_t\big]dW_t
\end{equation}
where $W_t$ designs a one-dimensional Brownian motion
\item A jump-equation
\begin{equation}\label{jump-equation}
d\rho_t=L(\rho_t)dt+\left[\frac{\mathcal{J}(\rho_t)}{Tr\big[\mathcal{J}(\rho_t)\big]}-\rho_t\right]
\left(d\tilde{N}_t-Tr\big[\mathcal{J}(\rho_t)\big]dt\right)
\end{equation}
where $\tilde{N}_t$ is a counting process with intensity $\int_0^t Tr\big[(\mathcal{J}(\rho_s)\big]ds$ (the different operators are described in the paper).
\end{itemize}

In this article, we shall focus on the case of jump-equation $(\ref{jump-equation})$, the diffusive case (\ref{difffff}) is treated in
details in \cite{Diffusion}. The equation (2) pose tedious problems in terms of mathematical justification. In the literature, the
question of existence and uniqueness of a solution is not really treated in details. Actually classical theorems cannot be applied
directly. Furthermore, the way of writing and presenting this equation is not clear. Indeed, how can we consider a driving process
which depends on the solution? With the expression (2), there is no intrinsic existence for such process. Even the notion of solution
is then not clear and a clearly probability space must be defined to give rigorous sense for such equations.

 Regarding the physical justification of Belavkin equation model, heuristic rules are usually used to derive these equations (see \cite{GST}).
   In order to obtain them rigorously, often a heavy
analytic machinery is used (Von-Neumann algebra, conditional expectation
in operator algebra, Fock space, quantum filtering...). This high
technology contrasts with the intuition of heuristics rules. In this article for the very first time,
 the continuous model (\ref{jump-equation})  is justified as a limit of a concrete discrete model of quantum repeated measurements.
 This approach, which has been already developed for the diffusive model (\ref{difffff}) in \cite{Diffusion}, is based on the model of
 quantum repeated interactions (see \cite{FRTC} for all details). The setup is the one of a small system interacting with an infinite
 chain of quantum system. Each piece of the chain interacts with the small system, one after the others, during a time $h$.
 After each interactions a measurement is performed on the piece of the chain which has interacted. The sequence of measurements induces then random perturbations of the small system. This is
  described by a Markov chain depending on the parameter $h$. Such Markov chain is called a \textit{discrete quantum trajectory}. In
   this article, by using a renormalization for the interaction (based on the result of \cite{FRTC}) and by considering the limit $h$
    goes to zero, we show that the quantum trajectory corresponding of the
     jump-equation can be obtained as continuous limit of particular discrete quantum trajectories.
 \bigskip

 As it is mentioned, similar results concerning existence, uniqueness and approximation
 for the diffusive equation is expressed in \cite{Diffusion}. The particularity of the jump-equation
 concerns the counting process $(\tilde{N}_t)$. In the diffusive case,
 the existence and uniqueness of solution concerns a \textit{strong} notion of solution. In the jump context,
 in a first time, we use the notion of \textit{Random Poisson Measure} to give a sense
 to the equation and in a second time, we deal with existence and uniqueness of a \textit{weak} solutions (i.e path by path).
 Concerning the convergence result in the diffusive case, the arguments are based on a convergence result for
 stochastic integrals due to Kurtz and Protter (\cite{WLT}, \cite{OWZ}). In the jump case, such techniques cannot be applied (essentially
 due to the stochastic intensity, see section 3). The result is then based on Random Coupling Limit
 method which allows to compare discrete quantum trajectories and continuous quantum trajectories.
 Finally to conclude to the convergence, we need the intermediary result of convergence of Euler scheme for the jump-equation. This result is also shown in this article because classical results cannot be applied in this context.
\bigskip

 This article is structured as follow.
we can
Section $1$ is devoted to remind the discrete model of repeated quantum measurements based on quantum repeated interactions models and
indirect measurement. Next, we make precise the model of a qubit in contact with a spin chain; this model corresponds to the effective
discrete model of Belavkin equations.

In Section 2, we investigate the study of continuous model of jump-equation. We present a way to define a concrete probabilistic
framework using random Poisson measure theory. Next, we solve problems of existence and uniqueness. The way, we describe the solution,
is next used in Section 3 to prove the convergence of discrete model.

In Section 3, we prove the convergence of a special discrete model to the jump continuous model.
 The arguments, we use in the jump case, are totally different of the diffusive case \cite{Diffusion}.
  Indeed the one, used in \cite{Diffusion}, are based on results of Kurtz and Protter (\cite{WLT}, \cite{OWZ}) and cannot
   be applied in this context. In order to prove the convergence result, we use explicit comparison by using random coupling method.
    Furthermore, we need to show another result of approximation of jump-equation: the Euler scheme convergence. This result is next use to prove the final convergence result.
\section{Quantum repeated measurements: A Markov chain}

\subsection{Quantum repeated measurements}

In this section, we present the mathematical model describing the quantum repeated measurements setup. Such model is treated in
details in \cite{Diffusion}. We do not remind all the details, we just present the rules and the Markov property of discrete quantum
trajectories.

The model is based on quantum repeated interactions model. As it is mentioned in the introduction, we consider a small system
$\mathcal{H}_0$ in contact with an infinite chain of quantum systems describing the environment. All the pieces of the chain are
identical and independent of each other; they are denoted by $\mathcal{H}$. One after the others, each copy of $\mathcal{H}$ interacts
with $\mathcal{H}_0$ during a time $h$. After each interaction, a measurement is performed on the copy $\mathcal{H}$, it involves then
a random perturbation of the state of $\mathcal{H}_0$.

Let start by describing a single interaction between $\mathcal{H}_0$ and $\mathcal{H}$ and the indirect measurement on $\mathcal{H}$.
Let $\mathcal{H}_0$ and $\mathcal{H}$ be finite dimensional Hilbert space. Each Hilbert space is endowed with a positive operator of
trace one. Such operators represent states of quantum systems. Let $\rho$ be the initial state on $\mathcal{H}_0$ and let $\beta$ be
the state of $\mathcal{H}$.

The coupling system is described by $\mathcal{H}_0\otimes\mathcal{H}$. The interaction is described by a total Hamiltonian $H_{tot}$ which is a self adjoint operator defined by
$$H_{tot}=H_0\otimes I+I\otimes H+H_{int}.$$
The operators $H_0$ and $H$ are the Hamiltonians of $\mathcal{H}_0$ and $\mathcal{H}$; they represent the free evolution of each
system. The operator
 $H_{int}$ is called the interaction Hamiltonian and describes the energy exchanges between the two systems. The Hamiltonian $H_{tot}$ gives rise to a unitary operator of evolution
$$U=e^{ih\,H_{tot}},$$
where $h$ is the time of interaction. Hence, after the interaction, in the \textit{Schr\"odinger picture}, the initial state
$\rho\otimes\beta$ on $\mathcal{H}_0\otimes\mathcal{H}$ becomes
$$\mu=U\,(\rho\otimes\beta)\,U^\star.$$
Let us now describe the measurement of an observable $A$ of $\mathcal{H}$. An observable is a self-adjoint operator and we consider
its spectral decomposition
$$A=\sum_{i=0}^p\lambda_i\,P_i.$$
Actually, we consider the observable $I\otimes A$ on $\mathcal{H}_0\otimes\mathcal{H}$. According to the law of quantum mechanics, the
measurement of $I\otimes A$ gives a random result concerning the eigenvalues of $I\otimes A$. It obeys to the following probability
law
$$P[\,\textrm{to observe}\,\,\lambda_i\,]=Tr\big[\,\mu\,I\otimes P_i\,\big].$$
After the measurement, if we have observed the eigenvalue $\lambda_i$ the reference state $\mu$ becomes
$$\mu_1(i)=\frac{I\otimes P_i\,\,\mu\,\,I\otimes P_i}{Tr\big[\,\mu\,I\otimes P_i\,\big]}.$$
Such phenomena is called \textit{Wave Packet Reduction Principle}. The state $\mu_1(i)$ is then the new reference state of
$\mathcal{H}_0\otimes\mathcal{H}$ conditionally to the observation of $\lambda_i$.

In general, one is only interested in the evolution of the small system $\mathcal{H}_0$. We use the partial trace operation to define
a state on $\mathcal{H}_0$ from a state on $\mathcal{H}_0\otimes\mathcal{H}$. The partial trace operation is defined as follows.

\begin{dthm}
Let $\alpha$ be a state on the tensor product $\mathcal{H}_0\otimes\mathcal{H}$. Hence there exists a unique state on $\mathcal{H}_0$,
denoted $\mathbf{E}_0[\alpha]$, which is characterized by the property
 $$\forall X \in \mathcal{B}(\mathcal{H}_0)\,\ Tr_{\mathcal{H}_0}[\,\eta
\,X\,]=Tr_{\mathcal{H}_0\otimes\mathcal{H}}[\,\alpha (X\otimes
I)\,],$$
where $Tr_{\mathcal{H}_0}$ corresponds to the trace of operator on $\mathcal{H}_0$ and similar definition
 for $Tr_{\mathcal{H}_0\otimes\mathcal{H}}.$
\end{dthm}
Now, let define the state $\rho_1(.)$ on $\mathcal{H}_0$ by
$$\rho_1(i)\,=\,\mathbf{E}_0[\mu_1(i)]\,,\,\,i=0,\ldots,p.$$
The state $\rho_1(.)$ is then the new reference state of $\mathcal{H}_0$ describing the result of one interaction and one measurement.
It represents a random state; each state $\rho_1(i)$ appears with probability $p_i=P[\textrm{to observe}\,\,\lambda_i]$. The
probability $p_i$ represents then the probability of transition from the state $\rho$ to the state $\rho_1(i)$.

As a consequence, $\mathcal{H}_0$ is now endowed with the state $\rho_1$ and a second copy of $\mathcal{H}$ can interact with
$\mathcal{H}_0$. In the same way, a measurement of $I\otimes A$ is then performed and by taking the partial trace, we get a new random
variable $\rho_2$ with similar transition probabilities (from $\rho_1$ to $\rho_2$) and so on. Hence we can define a random sequence
$(\rho_k)$ of states on $\mathcal{H}_0$; this sequence is called a \textit{discrete quantum trajectory}. It describes the random
evolution of the state of $\mathcal{H}_0$ undergoing quantum repeated interactions and quantum repeated measurements.

From the description of the rules of one interaction and one measurement and by construction of the random sequence $(\rho_k)$, the following proposition is straightforward.

\begin{pr}\label{RaS}
The discrete quantum trajectory $(\rho_k)$ is a Markov chain. The transition are described as follows. If $\rho_k=\chi_k$ then $\rho_{k+1}$ takes one of
the values
$$\mathbf{E}_0\left[\frac{I\otimes P_i\,U(\chi_k\otimes\beta)U^\star\,I\otimes
P_i}{Tr\big[\,U(\chi_k\otimes\beta)U^{\star}\,I\otimes P_i\,\big]}\right], \,\,i=0,\ldots,p,$$ with probability
$Tr\big[\,U(\chi_k\otimes\beta)U^\star\,P_i\,\big]$.
\end{pr}

A complete justification of this proposition and a construction of an appropriate probability space for the Markov chain can be found
in \cite{Diffusion}. All the details are not necessary in the following.

The next section is devoted to the study of a special case of a qubit in contact with a chain of spin.

\subsection{A qubit interacting with a chain of spin}

The mathematical setup describing such model is represented by $\mathcal{H}_0=\mathcal{H}=\mathbb{C}^2$. In this section, we make
explicit a discrete version of Belavkin equations.

Let us start by describing an evolution equation for the state $(\rho_k)$ in this context. In the case of $\mathbb{C}^2$, interesting
observable owns two different eigenvalues, that is $A=\lambda_0P_0+\lambda_1P_1$. The description of transitions of the Markov chain
$(\rho_k)$ can be expressed by the following equation
\begin{equation}\label{discreteeq}
\rho_{k+1}=\frac{\mathcal{L}_0(\rho_k)}{Tr\big[\mathcal{L}_0(\rho_k)\big]}\mathbf{1}^{k+1}_0+\frac{\mathcal{L}_1(\rho_k)}{Tr\big[
\mathcal{L}_1(\rho_k)\big]}\mathbf{1}^{k+1}_1,
\end{equation}
where for $i\in\{0,1\}$ terms $\mathcal{L}_i(\rho_k)$ corresponds to the "non-normalized" transition of $\rho_{k+1}$, that is
$\mathcal{L}_i(\rho_k)=\mathbf{E}_0\big[\,I\otimes P_i\,U(\rho_k\otimes\beta)U^\star\,I\otimes P_i\,\big]$. We denote
$p_{k+1}=Tr\big[\mathcal{L}_0(\rho_k)\big]=1-q_{k+1}$. Moreover, the term $\mathbf{1}^{k+1}_0$ corresponds to the random variable
which takes the value $1$ with probability $p_{k+1}$ and $0$ with probability $q_{k+1}$; it corresponds to the observation of
$\lambda_0$ at the k+1-th measurement. The inverse holds for the random variable $\mathbf{1}^{k+1}_1$.

In order to make more precise the evolution of the state $(\rho_k)$, we have to express in a more explicit way the terms
 $\mathcal{L}_i(\rho_k)$. For this, we introduce an appropriate basis. Let $(\Omega,X)$ denote an orthonormal basis of $\mathbb{C}^2$. For
$\mathcal{H}_{0}\otimes\mathcal{H}$,
  we consider the following basis $\Omega\otimes\Omega,X\otimes\Omega,\Omega\otimes
X,X\otimes X$. In this basis, the unitary-operator $U$ can be written as
$$U=\left( \begin{array}{cc}  L_{00} & L_{01} \\
  L_{10} & L_{11}
 \end{array}\right) $$
 where each $L_{ij}$ are operators on $\mathcal{H}_0$. For the state $\beta$, let choose the one dimensional projector on $\Omega$, that is
$$\beta=P_{\{\Omega\}}.$$
Furthermore if for $i\in\{0,1\}$, the projector $P_i=(p^i_{kl})_{k,l=0,1}$ in the basis $(\Omega,X)$, after computation we get
 \begin{eqnarray}\mathcal{L}_{i}(\rho_k)&=&\mathbf{E}_{0}\big[\,I\otimes
P_i\,\,U(\rho_k\otimes\beta)U^\star\,\,I\otimes P_i\,\big]\nonumber\\&=&p_{00}^i\,L_{00} \,\rho_k
\,L_{00}^\star+p_{01}^i\,L_{00}\,\rho_k \,L_{10}^\star+\,p_{10}^i\,L_{10}\,\rho_k \,L_{00}^\star+p_{11}^i\,L_{10}\,\rho_k
\,L_{10}^\star
\end{eqnarray}

In order to compare the discrete evolution and the continuous evolution in Section 3, we modify the way of writing
 the equation $(\ref{discreteeq})$ by introducing new random variables. For all $k\geq0$, let put $\nu_{k+1}=\mathbf{1}^{k+1}_1$ and define the random variable $$X_{k+1}=\frac{\nu_{k+1}-q_{k+1}}{\sqrt{q_{k+1}p_{k+1}}}.$$
We define the associated filtration on $\Sigma^{\mathbb{N}}$
$$\mathcal{F}_k=\sigma(X_i,i\leq k).$$ So by construction we have
$\mathbf{E}[X_{k+1}/\mathcal{F}_k]=0$ and $\mathbf{E}[X_{k+1}^2/\mathcal{F}_k]=1$. The random variables $(X_k)$ are then normalized
and centered. In terms of $(X_k)$, the discrete evolution equation for the discrete quantum trajectory becomes
\begin{equation}\label{discreteequation}
\rho_{k+1}=\mathcal{L}_{0}(\rho_k)+\mathcal{L}_{1}(\rho_k)+\left[-\sqrt{\frac{q_{k+1}}{p_{k+1}}}\mathcal{L}_{0}(\rho_k)
+\sqrt{\frac{p_{k+1}}{q_{k+1}}}\mathcal{L}_{1}(\rho_k)\right]X_{k+1}.
\end{equation}
This way, the evolution of the discrete quantum trajectory appears as a random perturbation of the equation
$\rho_{k+1}=\mathcal{L}_{0}(\rho_k)+\mathcal{L}_{1}(\rho_k)$ which describes actually the evolution without measurement. In
\cite{FRTC}, it is shown that $\mathcal{L}_0+\mathcal{L}_1$ is a completely
positive application and that it is a discrete Master equation. Equations of type (\ref{discreteequation}), obtained for different observable, can be then compared with Belavkin equations in terms of perturbation of Lindblad  Master equations (see
Section 2).

\textbf{Remark:} It is worth noticing that such equations depends only on the expression of eigen-projectors and do not depends on the
value of eigenvalues.

The following section is devoted to the study of the continuous Belavkin equations with jump.

\section{The jump Belavkin equation}

 First rigorous results and mathematical models describing evolution of systems undergoing quantum measurements are due to Davies in \cite{OQS}.
 From his works, one can derive heuristically the Belavkin equations. A heavy background is necessary in order to obtain these
  equations in a rigorous way (quantum filtering theory \cite{SSEMB},
instrumental process \cite{3MR2124562}). In this article, the continuous stochastic model is rigorously justified as a limit
of the discrete process given by the equation $(\ref{discreteequation})$. \\

In this article we focus on the jump-equation
\begin{equation}\label{jumpequation}
d\rho_t=L(\rho_t)dt+\left[\frac{\mathcal{J}(\rho_t)}{Tr\big[\mathcal{J}(\rho_t)\big]}-\rho_t\right]\left(d\tilde{N}_t-Tr\big[
\mathcal{J}(\rho_t)\big]dt\right)
\end{equation}
where $\tilde{N}_t$ is assumed to be a counting process with intensity $\int_0^t Tr\big[\mathcal{J}(\rho_s)\big]ds$.

Before to study this equation, let us speak briefly about the different operators appearing in the expression. In quantum physics, the operator $L$
is called the Lindbladian of the system. This is a classical generator of the dynamic of open quantum systems. It gives rise to the
Master equation (see the remark at the end of Section 1)
\begin{equation}\label{master}
\frac{d}{dt}\rho_t=L(\rho_t)=-i[H,\rho_t]-\frac{1}{2}\left\lbrace C^\star C,\rho_t\right\rbrace+C\rho_t C^\star
\end{equation}
 where $C$ is any $2\times 2$ matrix, the operator $H_0$ is the Hamiltonian of the qubit. In the equation $(\ref{jumpequation})$, the operator $\mathcal{J}$ is defined as
$\mathcal{J}(\rho)=C\,\rho\, C^\star$.

\subsection{The probability framework of the jump-equation}

The way of writing equation (\ref{jumpequation}) and defining the process $(\tilde{N}_t)$ is not absolutely rigorous in mathematical point of view. Indeed
in the definition, the process $(\tilde{N}_t)$ depends on existence of the solution of the jump-equation. However to treat the problem
of the existence of a solution, we have to consider firstly the definition of the driving process of the stochastic equation.

Actually, we cannot consider the existence of the process $(\tilde{N}_t)$ without the existence of the solution $(\rho_t)$ and
reciprocally. It imposes the following definition of solution for the jump-equation.

 \begin{df}\label{right}
Let $(\Omega,\mathcal{F},\mathcal{F}_t,P)$ be a filtered probabilistic space. A process-solution of $(\ref{jumpequation})$ is a
c\`adl\`ag process $(\rho_t)$ such that there exists a counting process $(\tilde{N}_t)$ with predictable compensator (or stochastic
intensity)
$$t\rightarrow\int_0^t Tr[\mathcal{J}(\rho_{s-})]ds$$ and such that
the couple $(\rho_t,\tilde{N}_t)$ satisfies almost surely
\begin{equation*}
\rho_t=\rho_0+\int_0^t\Big[L(\rho_{s-})-\mathcal{J}(\rho_{s-})+Tr\big[\mathcal{J}(\rho_{s-})\big]\,\rho_{s-}\Big]ds+\int_0^t
\left[\frac{\mathcal{J}(\rho_{s-})}{Tr[\mathcal{J}(\rho_{s-})]}-\rho_{s-}\right]d\tilde{N}_s.
\end{equation*}
\end{df}

This notion of solution comes from considerations of Jacod and Protter in
\cite{QR} \cite{CSPM}.

After this definition, the next step consists in constructing a process $\tilde{N}$ with stochastic intensity. For this, we use the
general theory of \textit{Random Measure} (for all details, see \cite{CSPM} or \cite{LTS}). Let us introduce this notion.

\begin{df}
Given a filtered probability space
$(\Omega,\mathcal{F},\mathcal{F}_{t},P)$, a random measure is a
family of measure $\mu=(\mu(\omega,.),\omega\in \Omega)$ on
$(\mathbf{R}_+\times
\mathbf{R}^d,\mathcal{B}(\mathbf{R}_{+})\otimes\mathcal{B}(\mathbf{R}^{d})).$

A random measure is said to be integer valued if\begin{enumerate}
\item For all $\omega\in\Omega$ $\mu(\omega,t\times\mathbf{R}^d)\leq
1$.
\item For all $A\in \mathcal{B}(\mathbf{R}_{+})\otimes\mathcal{B}(\mathbf{R}^{d})$,
the quantity $\mu(A)$ is valued in
$\mathbf{N}\bigcup\left\lbrace+\infty\right\rbrace$.
\end{enumerate}

\end{df}
\begin{df}A random Poisson measure on $(\Omega,\mathcal{F},\mathcal{F}_{t},P)$ is a
integer valued measure that verifies
\begin{enumerate}
\item The measure $m(A)=\mathbf{E}(\mu(A))$ on
$\mathcal{B}(\mathbf{R}_{+})\otimes\mathcal{B}(\mathbf{R}^{d})$  is
non atomic.
\item $m(0\times\mathbf{R}^{d})=0$.
\item If $t\in \mathbf{R}_{+}$ and if $A_i\in\mathcal{B}(\left] t,+\infty\right[)$,
$i=1,\ldots,l$ are two by two disjoint with $m(A_i)<+\infty$, the
random variables $\mu(A_i)$ are mutually independent and
independent from $ \mathcal{F}_t$.
\end{enumerate}
The measure $m$ is called the intensity of the random Poisson
measure $\mu$.
\end{df}

The following theorem shows how the random measure theory is used
to construct the process $(\tilde{N}_t)$.
\begin{thm}\label{process-solution}
Let $(\Omega,\mathcal{F},\mathcal{F}_{t},P)$ be a filtered probability space which supports a random Poisson measure $\mu$ on
$\mathbf{R}\times \mathbf{R}$ with intensity $dt\otimes dx$. Every process-solution of the equation
\begin{eqnarray}\label{Jump}
\rho_t&=&\rho_0+\int_0^t\Big[L(\rho_{s-})+Tr\big[\mathcal{J}(\rho_{s-})\big]\,\rho_{s-}-\mathcal{J}(\rho_{s-})\Big]ds\nonumber\\&&+\int_0^t\int_{\mathbf{R}_+}\left[\frac{\mathcal{J}(\rho_{s-})}{Tr[\mathcal{J}(\rho_{s-})]}-\rho_{s-}\right]\mathbf{1}_{0\leq
x\leq Tr[\mathcal{J}(\rho_{s-})]}\mu(ds, dx)
\end{eqnarray}
is a process-solution of equation $(\ref{jumpequation})$ satisfying Definition \ref{right}. For the process $(\tilde{N}_t)$, we put
\begin{equation}\label{pointprocess}
\tilde{N}_t=\int_0^t\int_{\mathbf{R}}\mathbf{1}_{0\leq x\leq Tr[\mathcal{J}(\rho_{s-})]}\mu(ds, dx).
\end{equation}
\end{thm}

 A general form of this theorem can be found in \cite{QR}.
  In this theorem, there are two parts. On the one hand, we must prove that the process
given by $(\ref{pointprocess})$ is well defined, that is, it is a non-explosive process. On the other hand, we must prove that any
solution of equation $(\ref{Jump})$ satisfies Definition \ref{right}.

The non-explosive property of $(\tilde{N}_t)$ is related to the boundness character of the stochastic intensity $t\rightarrow
Tr[\mathcal{J}(\rho_{t-})]$. Here, a straightforward computation shows that there exists a constant $K$ such that for all state $\rho$,
we have $0\leq Tr[\mathcal{J}(\rho)]\leq K$. It implies directly that for all t, the quantity
$Tr[\mathcal{J}(\rho_{t-})]=\lim_{s<t,s\rightarrow t} Tr[\mathcal{J}(\rho_s)]$ satisfies
$$0\leq Tr[\mathcal{J}(\rho_{t-})]\leq K,
$$
 the intensity is then bounded. This property appears also in the proof of the theorem.\\
\begin{pf}
Let show that the counting process $(\tilde{N}_t)$ given by $(\ref{pointprocess})$ is non explosive. For all c\`adl\`ag matricial
process $(X_t)$, we define the time explosion
$$T^X=\inf\left(t:\tilde{N}^X_t=+\infty\right).$$
Let show that, if $(\rho_t)$ takes values in the set of states, the explosion time $T^\rho=\infty$ almost surely. For this, we
introduce the following stopping times
$$T_n=\inf\left(t, \int_0^tTr[\mathcal{J}(\rho_{s-})]ds\geq n\right),\,\,n\geq1.$$
It was clear that $\int_0^{T_n}Tr[\mathcal{J}(\rho_{s-})]ds\leq n$ and thanks to the property of Poisson random measure
\begin{eqnarray*}\mathbf{E}\left[\tilde{N}^\rho_{T_n}\right]&=
&\mathbf{E}\left[\int_0^{T_n}\int_{\mathbf{R}}\mathbf{1}_{0\leq x\leq
Tr[\mathcal{J}(\rho_{s-})]}dsdx\right]\\&=&\mathbf{E}\left[\int_0^{T_{N}}Tr\left[\mathcal{J}(\rho_{s-})\right]ds\right].
\end{eqnarray*}
As a consequence we have $T_n\leq T^\rho$ almost surely. It is then sufficient to show that $\lim T_n=\infty$. The property $0\leq
Tr[\mathcal{J}(\rho_t)]\leq K$ implies $\int_0^tTr[\mathcal{J}(\rho_{s-})]ds<\infty$, hence we have
$\lim_{n\rightarrow\infty}T_n=+\infty$. Finally we have construct a counting process without explosion for all c\`adl\`ag process
which takes values in the set of states. Concerning the property that a process-solution of $(\ref{Jump})$ satisfies Definition
\ref{right}, this result follows from the construction of $(\tilde{N}_t)$.
\end{pf}\\
\hphantom{c}\\
\textbf{Remark:} This remark concerns a way to construct a random Poisson measure with uniform
intensity measure. Let
$(\Omega,\mathcal{F},P)$ be the probability space of a Poisson point process $N$ on
$\mathbf{R}\times\mathbf{R}$. The natural random Poisson measure attached
with $N$ is defined for all
$A\in\mathcal{B}(\mathbf{R})\otimes\mathcal{B}(\mathbf{R})$ by
$$\mu(.,A)=N(.,A).$$
For all Borel subset $A\in\mathcal{B}(\mathbf{R})\otimes\mathcal{B}(\mathbf{R})$, we have $\mathbf{E}[N(.,A)]=\lambda(A)$ where
$\lambda$ denotes the Lebesgue measure. The Poisson random measure $\mu$ satisfies then the conditions of Theorem
$\ref{process-solution}$. This particular random Poisson measure is used in Section $3$ to realize continuous quantum trajectories and
discrete quantum trajectories in a same probability space and to prove the convergence result.

After defining this clear probability framework, it remains to deal with problems of existence and uniqueness of solution.

\subsection{Existence and uniqueness}

Let $(\Omega,\mathcal{F},\mathcal{F}_t,P)$ which supports a random Poisson measure $\mu$. In order to solve problems of existence and
uniqueness of solution for equations of type $(\ref{Jump})$, a way consists firstly in solving the ordinary differential part
\begin{equation}\label{Ode}
\rho_t=\rho_0+\int_0^t\Big[L(\rho_{s-})+Tr\big[\mathcal{J}(\rho_{s-})\big]\,\rho_{s-}-\mathcal{J}(\rho_{s-})\Big]ds,
\end{equation}
and secondly defining times of jump.  Classical theorems are usually based on Lipschitz property (cf \cite{PLL},\cite{QR} or \cite{OTS}). However in our case, the ordinary differential equation part contains non-Lipschitz functions. Hence the application of classical theorems is not straightforward.

The first step is then solving a non-Lipschitz ordinary differential equation (\ref{Ode}). For this, we use the fact that
this equation preserves the property of being a \textit{pure state}. Remind that pure states are a particular class of states which
are one dimensional projectors. The property of preserving pure states for $(\ref{Ode})$ is expressed as follows.

 \begin{pr}\label{equivalence}
 Let x be any vector of norm one in $\mathbb{C}^2$, if the Cauchy problem
 \begin{equation}\label{cauchy}\left\{\begin{array}{ccl}
dx_t&=&\left[-iH_0-\frac{1}{2}C^\star
C+\frac{1}{2}\eta_t\right]x_tdt \\
  x_0&=&x\end{array}\right.\end{equation}
  where $\eta_t=\langle x_t,C^\star Cx_t\rangle$
has a solution then $\Vert x_t\Vert=1$ for all $t>0$.

Furthermore the process $(\rho_t)$ of one-dimensional projector defined by $\rho_t=P_{\{x_t\}}$ for all $t>0$ is solution of Cauchy
problem
   \begin{equation}\left\{\begin{array}{ccl}
d\rho_t&=&\Big[L(\rho_{t})+Tr\big[\mathcal{J}(\rho_{t})\big]\rho_{t}-\mathcal{J}(\rho_{t})\Big]dt \\
  \rho_0&=&P_{\{x\}}.\end{array}\right.\end{equation}
 \end{pr}
\begin{pf}
 Let $(x_t)$ be the solution of $(\ref{cauchy})$, thanks to the fact that $H_0$
is self-adjoint and $\eta_t=\langle x_t,C^\star Cx_t\rangle=\langle Cx_t, Cx_t\rangle$, a straightforward computation gives
 \begin{eqnarray*}
 \frac{d}{dt}\langle x_t,x_t\rangle&=&\langle\frac{d}{dt}x_t,x_t\rangle+\langle
x_t,\frac{d}{dt}x_t\rangle\\
 &=&\langle\left[-iH_0-\frac{1}{2}C^\star
C+\frac{1}{2}\eta_t\right]x_t,x_t\rangle+\\&&\langle x_t,\left[-iH_0-\frac{1}{2}C^\star
C+\frac{1}{2}\eta_t\right]x_t\rangle\\
 &=&-\frac{1}{2}\eta_t+\frac{1}{2}\eta_t\langle x_t,x_t\rangle.
 \end{eqnarray*}As a consequence, if $\langle
x_0,x_0\rangle=1$, then $\langle x_t,x_t\rangle=1$ for all $t$. Let $\rho_t=P_{\{x_t\}}$, for all $y$, we have
$$\rho_t\,y=\langle x_t,y\rangle x_t.$$ We can derive $d\rho_ty$ and we have
\begin{eqnarray*}
\frac{d}{dt}\,\rho_ty&=&\langle\frac{d}{dt}\,x_t,y\rangle
x_t+\langle
x_t,y\rangle\frac{d}{dt}\,x_t\\
&=&\langle \left[-iH_0-\frac{1}{2}C^\star C+\frac{1}{2}\eta_t\right]x_t,y\rangle x_t+\langle x_t,y\rangle\left[-iH_0-\frac{1}{2}C^\star
C+\frac{1}{2}\eta_t\right]x_t\\
&=&\langle x_t,\left[-iH_0-\frac{1}{2}C^\star C+\frac{1}{2}\eta_t\right]^\star y\rangle x_t+\left[-iH_0-\frac{1}{2}C^\star
C+\frac{1}{2}\eta_t\right]\langle x_t,y\rangle x_t\\&=&\rho_t\left[-iH_0-\frac{1}{2}C^\star C+\frac{1}{2}\eta_t\right]^\star
y+\left[-iH_0-\frac{1}{2}C^\star
C+\frac{1}{2}\eta_t\right]\rho_t y\\
&=&\Big[L(\rho_{t})+Tr\big[\mathcal{J}(\rho_{t})\big]\,\rho_{t}-\mathcal{J}(\rho_{t})\Big]y
\end{eqnarray*}
and the result follows.\end{pf}
\bigskip

With this expression in term of pure states, we can show that equation $(\ref{Ode})$ admits a solution. For this, we need the
following characterization of pure states in $\mathbb{C}^2$.
 \begin{lem}\label{onedimensionnal}
 Let $\rho$ be a state on $\mathbf{C}^2$. If there exists a vector $x\in\mathbf{C}^2$
such that $\langle x;\rho x\rangle=0$, the state $\rho$ is a one dimensional projector.
 \end{lem}
We do not give the proof, it is a simple linear algebra result.

The following proposition expresses the result concerning the existence of a solution for $(\ref{Ode})$.

 \begin{pr}\label{Cauchyproblem}
 Let $\rho$ be any state, the Cauchy problem
   \begin{equation}\label{Cauchy}\left\{\begin{array}{ccl}
d\rho_t&=&\Big[L(\rho_{t})+Tr\big[\mathcal{J}(\rho_{t})\big]\rho_{t}-\mathcal{J}(\rho_{t})\Big]dt \\
  \rho_0&=&\rho\end{array}\right.\end{equation}
   has a unique solution defined for all time $t$.

  Furthermore, if there exists $t_0$ such that $\rho_{t_0}$ is a one dimensional
projector, the solution of $(\ref{Cauchy})$ after $t_0$ is valued in the set of pure states.
 \end{pr}
\begin{pf}
 As the coefficients are not Lipschitz, the Theorem of Cauchy
Lipschitz cannot be applied directly. However the coefficients are $C^{\infty}$, so locally
Lipschitz and we can use a truncation method. The ordinary
equation is of the following form
\begin{equation*}d\rho_t=f(\rho_t)\,dt
\end{equation*}
where $f$ is $C^\infty$ and
$f(A)=L(A)+Tr[\mathcal{J}(A)]A-\mathcal{J}(A)$. We define the
truncation function $\varphi$ from $\mathbf{R}$ to $\mathbf{R}$
defined by
 $$\varphi_k(x)=\left\{\begin{array}{cc} -k &\textrm{if $x\leq-k$}\\ x & \textrm{if
$-k\leq x\leq k$}\\ k & \textrm{if $-k\leq x\leq k$}.\end{array}\right.$$ For a matrix $A=(a_{ij})$, we define by extension
$\tilde{\varphi}_k(A)=\varphi_k(Re(a_{ij}))+i\varphi_k(Im(a_{ij}))$. Thus $f\circ\tilde{\varphi}_k$ is Lipschitz. Now we consider the
truncated equation
$$d\rho_{k,t}=f\circ\tilde{\varphi}_k(\rho_{k,t})dt.$$
The Cauchy-Lipschitz Theorem can be applied because $f\circ\tilde{\varphi}_k$ is Lipschitz and there exists a unique solution
$t\mapsto\rho_{k,t}$ defined for all t.

Now, we can define the following time
$$T_k=inf\{t,\exists(ij)/\,\vert(Re(a_{ij}(\rho_{k,t})))\vert=k\,\,or\,\,\vert(Im(a_{ij}(\rho_{k,t})))\vert=k\}.$$
As $\rho_0$ is a state, if k is chosen large enough we have $T_k>0$
  and for all $t\leq T_k$, we have $\tilde{\varphi}_k(\rho_{k,t})=\rho_{k,t}$.
   Thus $t\mapsto\rho_{k,t}$ is the unique solution of the ordinary equation
$(\ref{Ode})$ (without truncation) on $[0,T_k]$. The classical method in order to solve an equation with non Lipschitz
coefficients is to put $T=\lim_k T_k$
     and to show that $T=\infty$.
     
      Here the situation is more simply because if $\rho_0$
      is a state, we can show that the solution is valued on the set of
states. Indeed, as $\Vert \rho\Vert\leq1$ when $\rho$ is a state, we have for example $\tilde{\varphi}_2(\rho)=\rho$. Let show that
the solution obtained by truncation is valued in the set of states, it will imply that $T_2=\infty$. We have to show that the solution
is self adjoint, positive and of trace one.

Let start with the self adjoint property and trace property. On $[0,T_2]$, as the ordinary differential equation is Lipschitz, we can
solve it by Picard method. Let define
           \begin{equation}\left\{\begin{array}{ccl}\rho_{n+1}(t)&=&\rho_n(0)+
     \int_0^tf\circ \tilde{\varphi}_k(\rho_n(s))ds\\
     \rho_0(t)&=&\rho.
     \end{array}\right.
     \end{equation}
 It is easy to show with the right definition of $f$ that this sequence
             is self adjoint with trace one. Hence it proves (by convergence of the Picard method) that for all $t\leq T_2$, the matrix
$\rho_{2,t}$ is self adjoint with trace one.

 It remains to prove the positivity property which cannot actually be proved by the previous Picard method approach.
 This condition is nevertheless a consequence
of Proposition $\ref{equivalence}$ and Lemma $\ref{onedimensionnal}$. We must prove that for all $y\in\mathbf{R}^{2}$ and for all
$t\leq T_2$, we have $\langle y,\rho_{2,t}y\rangle\geq0$. Let define $$T^0=\inf\{t\leq T_2,\exists y\in\mathbf{R}^{2}/\langle
y,\rho_{2,t}y\rangle=0\}.$$

If $T^0=T_2$, an argument of continuity for solution of ordinary differential equation implies that $\langle y,\rho_{2,t}y\rangle\geq
0$ for all $t\leq T_2$
         and all $y\in\mathbf{C}^{2}$. The solution is then valued on the set of states.

If $T^0<T_2$,  by continuity there exists
          $x\in\mathbb{R}^{2}$  such that $\langle x,\rho_{2,T^0}x\rangle=0$ and for all $t\leq T^0$ and for all
           $y\in\mathbf{R}^{2}$, we have $\langle y,\rho_{2,t}y\rangle\geq0$.
            It means that on $[0,T^0]$, the solution $t\mapsto \rho_{2,t}$ is valued
on the set of states. Moreover for some $x$, we have $\langle x,\rho_{2,T^0}^x\rangle=0$.
             Thanks to Lemma $\ref{onedimensionnal}$, the operator $\rho_{2,T^0}$ is a one dimensional projector. Hence, we can now consider the ordinary differential equation with initial state $\rho_{T^0}=\rho_0$. We are face to the Cauchy problem
$(\ref{cauchy})$ which is equivalent to the problem $(\ref{Cauchy})$ thanks to Proposition $\ref{equivalence}$. This problem can be
also solved by truncation method, the fact that the norm is conserved implies that the solution is defined for all $t$ (the truncation
is actually not necessary). Thanks to Proposition $\ref{equivalence}$, we have a solution for the initial Cauchy problem
$(\ref{Cauchy})$ which defines a state-process. We have then finally proved that on $[T^0,T_2]$ the solution is valued on the set of
states.

On $T_2$, the solution is then a state. A local argument and the uniqueness in the Cauchy-Lipschitz theorem allows us to conclude that $T_2=\infty$ and that there
exists a unique solution of the ordinary differential equation $(\ref{Ode})$.
 \end{pf}
 \bigskip

 This above
proposition is essential in the proof of the following theorem concerning existence and uniqueness of the jump-equation $(\ref{Jump})$.

 \begin{thm} Let $(\Omega,\mathcal{F},\mathcal{F}_t,P)$ be a probability space which supports
 a Poisson random measure $\mu$ whose the intensity measure is
$dx\otimes dt$. Let $\rho_0$ be any state, the jump-Belavkin
equation
 \begin{eqnarray}\label{Good}
\rho_t&=&\rho_0+\int_0^t\Big[L(\rho_{s-})+Tr\big[\mathcal{J}(\rho_{s-})\big]\rho_{s-}-\mathcal{J}(\rho_{s-})\Big]ds\nonumber\\&&+\int_0^t\int_{\mathbf{R}_+}\left[\frac{\mathcal{J}(\rho_{s-})}{Tr[\mathcal{J}(\rho_{s-})]}-\rho_{s-}\right]\mathbf{1}_{0\leq
x\leq Tr[\mathcal{J}(\rho_{s-})]}\mu(ds, dx)
\end{eqnarray}
 admits a unique solution defined for all
time. The process-solution $(\rho_t)$ takes values on the set of states.
 \end{thm}
 \begin{pf} Such equation is solved paths by paths. As the initial condition $\rho_0$ is a state, Proposition $\ref{Cauchyproblem}$ states the existence and the uniqueness for the Cauchy problem $(\ref{Cauchy})$. Let define the
first time of jump, we put
\begin{equation}\left\{\begin{array}{ccl}
\rho(1)_t&=&\rho_0+\int_0^tf(\rho(1)_s)ds \\
T_1&=&\inf\{t,\tilde{N}^{\rho(1)}_t>0\}.\end{array}\right.\end{equation} The first equality expresses the fact that the Cauchy problem
admits a solution. The second equality defines a random variable $T_1$ which is the stopping time of the first jump. We have
\begin{equation*}\tilde{N}^{\rho(1)}_{T_1}(\omega)=\mu(\omega,G(\rho,T_1,0))=1,\end{equation*}
where $G(\rho,t,s)=\{(u,y)\in\mathbf{R}^2/t<u<s,0<y<Tr[\mathcal{J}(\rho_u)]\}$. The quantity $\mu(\omega,G(\rho,t,s))$ represents the
number of point under the curves $t\rightarrow Tr[\mathcal{J}(\rho_t)]$.

Now if $T_1=\infty$, the solution of the jump-equation is given by the solution of the ordinary differential equation; there are no
jumps.

If $T_1<\infty$, let construct the second jump time and explain how the solution of the jump-equation is defined. On $[0,T_1[$, we put
$\rho_t=\rho(1)_t$ and at the jump time $T_1$, we implement the value of the jump, that is,
$$\rho_{T_1}=\rho_{T_1-}+\left[\frac{\mathcal{J}(\rho_{T_1-})}{Tr[\mathcal{J}(\rho_{T_1-})]}-\rho_{T_1-}\right]=\frac{\mathcal{J}(\rho_{T_1-})}{Tr[\mathcal{J}(\rho_{T_1-})]}.$$
This matrix is well defined because the property $T_1<\infty$ implies $Tr[\mathcal{J}(\rho_{T_1-})]>0$. Moreover, the matrix
$\rho_{T_1}$ is a state. The ordinary differential equation $(\ref{Ode})$ can be solved again with this initial condition. We can then
define
  \begin{equation}\left\{\begin{array}{ccl} \rho(2)_t&=&\rho(1)_t\,\,
\textrm{on}\,\, [0,T_1[\\\rho(2)_t&=&\rho_{T_1}+\int_{T_1}^tf(\rho(2)_s)ds \\
T_2&=&\inf\{t>T_1,\tilde{N}^{\rho(2)}_t>\tilde{N}^{\rho(1)}_{T_1}\}.\end{array}\right.\end{equation}
  The random variable $T_2$ is then the second jump-time. If
$T_2<\infty$, we have $Tr[\mathcal{J}(\rho(2)_{T_2-})]>0$ and we continue in the same way. The solution of the jump-equation is then defined as the solution of the ordinary differential
equation between the jump times and at jump times we put the value of jump.

More precisely, we define a sequence of random times $T_n$ and a sequence of processes
  \begin{equation}\left\{\begin{array}{ccl} \rho(n)_t&=&\rho(n-1)_t\,\,
\textrm{on}\,\, [0,T_{n-1}[\\\rho_{T_{n-1}}&=&\frac{\mathcal{J}(\rho_{T_{n-1}-})}{Tr[\mathcal{J}(\rho_{T_{n-1}-})]}
\\\rho(n)_t&=&\rho_{T_{n-1}}+\int_{T_{n-1}}^tf(\rho(n-1)_s)ds \\
T_n&=&\inf\{t>T_{n-1},\tilde{N}^{\rho(n)}_t>\tilde{N}^{\rho(n-1)}_{T_1}\}.\end{array}\right.\end{equation}
  All the processes are well defined because the matrices defined at each jump time are states and the Cauchy problem can be solved. The sequence of random
stopping times $(T_{n})$ satisfies $T_{n+1}>T_n$ on the set
$\{T_n<\infty\}$. In a classical way, we put
$$T=\lim_{n\rightarrow\infty}T_n,$$
and we can define the process-solution $(\rho_t)$ of the jump-Belavkin equation on $[0,T[$. For all $t<T$, we put
\begin{equation}\rho_t=\rho(n)_t \,\,\textrm{on}\,\, [0,T_n[.
\end{equation}
This process is clearly a solution of the jump-Belavkin equation $(\ref{Good})$ and it is valued on the set of states. The uniqueness
is implied by the uniqueness of the solution of Cauchy problem (cf Proposition $\ref{Cauchyproblem}$). Moreover any other solution
is forced to have the same random jump-time, it implies the uniqueness.

 In order to finish the proof, we must show $T=\infty$ a.s. This random time is the
explosion-time of $\tilde{N}^\rho$, we cannot directly apply the
result of the theorem $3$ because the definition for all $t$ of
$(\rho_t)$ is not already proved. However $(\rho_t)$ is a state
valued process, so we have $Tr[\mathcal{J}(\rho_t)]<K$ for all $t$
then 
$$\mathbf{E}\left[\tilde{N}^\rho_{T_p\wedge
n}\right]\leq\mathbf{E}\left[\tilde{N}^\rho_ n\right]\leq Kn.$$ Furthermore $\tilde{N}^\rho_{T_p\wedge n}=p$ on $\{T_P<n\}$, it
follows that $pP[T_p<n]\leq Kn$, then we have $P[T\leq n]=0$ for all $n$ and the result is proved.
\end{pf}\\

The proof of theorem gives an explicit way to construct the solution and the random times of jumps. This construction can be used in a
more general context. In the literature, the problem of existence is sometimes taken into account by a way using filtering method. A
starting point is using a linear equation and next a change of measure. However this method do not give an explicit construction of
the solution whereas we need it to prove the final convergence result.
\bigskip

In order to prove the final result, we use random coupling method. For this, we need an explicit realization of the process
$(\rho_t)$. For this, we use Poisson point process theory. Let consider the probability space $(\Omega,\mathcal{F},\mathcal{F}_t,P)$ of
a Poisson point process $N$ on $\mathbf{R}\times\mathbf{R}$. Let denote $N(\omega,ds,dx)$, the differential increment of the random
measure $N$. The continuous quantum trajectory $(\rho_t)$ satisfies
 \begin{eqnarray}\label{integral}\rho_t&=&\rho_0+\int_0^t\Big[L(\rho_{s-})+Tr\big[\mathcal{J}(\rho_{s-})\big]\rho_{s-}-\mathcal{J}(\rho_{s-})\Big]ds\nonumber\\&&+\int_0^t\int_{[0,K]}\left[\frac{\mathcal{J}(\rho_{s-})}{Tr[\mathcal{J}(\rho_{s-})]}-\rho_{s-}\right]\mathbf{1}_{0\leq
x\leq Tr[\mathcal{J}(\rho_{s-})]}N(.,ds, dx).
\end{eqnarray}
It is interesting to notice that we can work on $[0,K]$ because $Tr[\mathcal{J}(\rho_t)]\leq K$ for all process valued in the set of
states. Based on this property, the following remark gives another equivalent way to define the solution of $(\ref{Jump})$.\\
\textbf{Remark:} The function $t\rightarrow
card(N(.,[0,K]\times[0,t]))=\mathcal{N}_t$ defines a standard
Poisson process with intensity $K$. Thus for the filtration
$\mathcal{F}_t$, we can choose the natural filtration of this
process. The Poisson random measure and the previous process
generate on $[0,T]$ (for a fixed $T$) a sequence
$\{(\tau_i,\xi_i),i\in\{1,\ldots,\mathcal{N}_t)\}\}$ where each
$\tau_i$ represents one jump time of $\mathcal{N}$. Moreover, the
random variables $\xi_i$ are random uniform variables on $[0,K]$.
Consequently, we can write our continuous quantum trajectory in the
following way
 \begin{eqnarray*}\rho_t&=&\rho_0+\int_0^t\Big[L(\rho_{s-})+Tr\big[\mathcal{J}(\rho_{s-})\big]\rho_{s-}-\mathcal{J}(\rho_{s-})\Big]ds\\&&+\sum_{i=1}^{\mathcal{N}_t}\left[\frac{\mathcal{J}(\rho_{\tau_i-})}{Tr[\mathcal{J}(\rho_{\tau_i-})]}-\rho_{\tau_i-}\right]\mathbf{1}_{0\leq
\xi_i\leq Tr[\mathcal{J}(\rho_{\tau_i-})]}.
\end{eqnarray*}

The next section concerns the final convergence result.
\section{Approximation and convergence theorems}

This section is devoted to convergence theorems. The discrete quantum trajectory defined in Section $1$ is shown to converge to the
continuous quantum trajectory which is solution of the jump Belavkin equation $(\ref{integral})$.

The way to proceed to show this result is the following. Let $h=1/n$ be the time of interaction between $\mathcal{H}_0$ and a copy of
$\mathcal{H}$. This way, we shall define discrete quantum trajectories depending on parameter $n$. Next, based on result of Attal
Pautrat in \cite{FRTC}, we introduce asymptotic assumptions on the unitary operator $U$ describing the evolution. Next, for special
observable, we show that the corresponding discrete quantum trajectory is a  "good candidate" to converge to the solution of the
jump-equation. To prove the final result, it needs two step.
\begin{itemize}
\item The first step is to realize the discrete quantum trajectory (in asymptotic form) in the same probability space of the continuous quantum trajectory.
\item The second step is to compare the discrete quantum trajectory with an Euler scheme of the jump-equation.
\end{itemize}
\textbf{Remark:} In the literature, there exist classical Euler scheme results for jump-equation (see \cite{OTS}). At our knowledge, the
jump-Belavkin equation do not come into these classical setups. As a consequence, we establish the convergence of Euler scheme in our
context.

Let start by investigating the asymptotic behavior of discrete quantum trajectories.

\subsection{Asymptotic form of discrete jump-Belavkin process}

In this section, we present the asymptotic assumption concerning the interaction. Next we translate this assumption in equations
\begin{equation} \rho_{k+1}=\mathcal{L}_{0}(\rho_k)+\mathcal{L}_{1}(\rho_k)+\left[
-\sqrt{\frac{q_{k+1}}{p_{k+1}}}\mathcal{L}_{0}(\rho_k)+\sqrt{\frac{p_{k+1}}{q_{k+1}}}\mathcal{L}_{1}(\rho_k)\right] X_{k+1}.
\end{equation}

Consider a partition of $[0,T]$ in subintervals of equal size
$1/n$. The unitary evolution depends naturally on $n$ and we put
$$U(n)=\left( \begin{array}{cc}  L_{00}(n) & L_{01}(n) \\
  L_{10}(n) & L_{11}(n)
 \end{array}\right).$$
 In the quantum repeated interactions setup \cite{FRTC}, it is shown that the coefficients $L_{ij}$ must obey precise asymptotic to
 obtain non-trivial limit when $n$ goes to infinity. More precisely, they have shown that quantum stochastic differential equations (also called \textit{"Hudson-Parthasarathy equations"}),
 describing continuous time interaction models, can be obtained as continuous limit models of quantum repeated interactions by rescaling
 discrete interactions (such results use Toy Fock space and Fock space formalism).

We do not need the total result of \cite{FRTC}, we just need the expression of $L_{ij}$. In our context the asymptotic are the following
\begin{eqnarray}L_{00}(n)&=&I+\frac{1}{n}(-iH-\frac{1}{2}CC^\star)+\circ\left(\frac{1}{n}\right)\\
  L_{10}(n)&=&\frac{1}{\sqrt{n}}C+\circ\left(\frac{1}{n}\right).
  \end{eqnarray}
  Concerning the description of the total Hamiltonian $H_{tot}(n)$, we can write  \begin{eqnarray*}H_{tot}(n)&=&H_0\otimes I+I\otimes\left(\begin{array}{cc} 1 &
0\\0&0\end{array}\right)+\frac{1}{\sqrt{n}}\left[C\otimes\left(\begin{array}{cc}
0 & 0\\1&0\end{array}\right)+C^\star\otimes\left(\begin{array}{cc}
0 & 1\\0&0\end{array}\right)\right]+\circ\left(\frac{1}{n}\right).\end{eqnarray*}
With the time discretization, we obtain a discrete process depending on $n$
  \begin{eqnarray}\label{assy}
\rho_{k+1}(n)&=&\mathcal{L}_{0}(n)(\rho_k(n))+\mathcal{L}_{1}(n)(\rho_k(n))\nonumber\\&&
+\left[
-\sqrt{\frac{q_{k+1}(n)}{p_{k+1}(n)}}\mathcal{L}_{0}(n)(\rho_k(n))+\sqrt{\frac{p_{k+1}(n)}{q_{k+1}(n)}}\mathcal{L}_{1}(n)(\rho_k(n))\right]
X_{k+1}(n).
\end{eqnarray}
Remind that by definition, for $(X_k)$, we have
\begin{equation}X_{k+1}(n)(i)=\left\{\begin{array}{ccl}
-\sqrt{\frac{q_{k+1}(n)}{p_{k+1}(n)}}&\textrm{with
probability}&p_{k+1}(n)\,\,\textrm{if}\,\, i=0\\
\sqrt{\frac{p_{k+1}(n)}{q_{k+1}(n)}}&\textrm{with
probability}&q_{k+1}(n)\,\,\textrm{if}\,\, i=1\end{array}\right.,
\end{equation}
where $p_{k+1}(n)=Tr\Big[I\otimes P_0\,\,U(n)(\rho_k\otimes\beta)U(n)^\star I\otimes P_0\Big]=1-q_{k+1}(n)$ (remind that $P_0$ is one
of the eigen-projector of the measured observable). Hence, the continuous limit behavior of the sequence $(\rho_k)$ will depend on
limit behavior of the random variables $(X_k)$ and then depends on asymptotic behavior of probabilities $p_{k+1}$ and $q_{k+1}$.

By applying the asymptotic assumption in equation $(\ref{assy})$, depending on the observable, we get two different situations.
\begin{itemize}
\item If the observable $A$ is diagonal in the orthonormal basis $(\Omega,X)$ of $\mathcal{H}$, that is $A=\lambda_0\left( \begin{array}{cc}  1 & 0 \\
  0 & 0
 \end{array}\right)+\lambda_1\left( \begin{array}{cc}  0 & 0 \\
  0 & 1
 \end{array}\right)$. We obtain the asymptotic for the probabilities
 \begin{eqnarray*}
p_{k+1}(n)&=&1-\frac{1}{n}Tr[\mathcal{J}(\rho_k(n))]+\circ\left(\frac{1}{n}\right)\\
q_{k+1}(n)&=&\frac{1}{n}Tr[\mathcal{J}(\rho_k(n))]+\circ\left(\frac{1}{n}\right).
\end{eqnarray*}
The discrete equation becomes
\begin{eqnarray*}
\rho_{k+1}(n)-\rho_k(n)&=&\frac{1}{n}L(\rho_k(n))+\circ\left(\frac{1}{n}\right)\\
&&+\left[\frac{\mathcal{J}(\rho_k(n))}{Tr(\mathcal{J}(\rho_k(n)))}-\rho_k(n)+\circ(1)\right] \sqrt{q_{k+1}(n)p_{k+1}(n)} \,X_{k+1}(n).
\end{eqnarray*}
\item If the observable is non diagonal in the basis $(\Omega,X)$, we consider $P_0=\left(
\begin{array}{cc}
p_{00} & p_{01} \\
  p_{10} & p_{11}
 \end{array}\right)$ and $P_1=\left( \begin{array}{cc}  q_{00} & q_{01} \\
  q_{10} & q_{11}
 \end{array}\right)$ we have
 \begin{eqnarray*}
p_{k+1}&=&p_{00}+\frac{1}{\sqrt{n}}Tr[\rho_k(p_{01}C+p_{10}C^\star)]+\frac{1}{n}Tr[\rho_k(p_{00}(C+C^\star))]+\circ\left(\frac{1}{n}\right)\\
q_{k+1}&=&q_{00}+\frac{1}{\sqrt{n}}Tr[\rho_k(q_{01}C+q_{10}C^\star)]+\frac{1}{n}Tr[\rho_k(q_{00}(C+C^\star))]+\circ\left(\frac{1}{n}\right).
\end{eqnarray*}
The discrete equation becomes then
\begin{eqnarray*}
\rho_{k+1}-\rho_k&=&\frac{1}{n}L(\rho_k)+\circ\left(\frac{1}{n}\right)+
\left[e^{i\theta}C\rho_k+e^{-i\theta}\rho_kC^\star\right.\\&&\left.-Tr[\rho_k(e^{i\theta}C+e^{-i\theta}C^\star)]\,\rho_k+\circ
(1)\right]\frac{1}{\sqrt{n}}X_{k+1}.
\end{eqnarray*}
\end{itemize}

From these descriptions, we can define processes $\rho_{[nt]}(n)$ by
\begin{eqnarray}
\rho_{[nt]}(n)&=&\rho_0+\sum_{k=1}^{[nt]-1}\rho_{k+1}-\rho_k\nonumber\\
&=&\rho_0+\sum_{k=1}^{[nt]-1}\bigg(L(\rho_k)+\circ(1)\Big)\frac{1}{n}+\sum_{k=1}^{[nt]-1}\mathcal{Q}_i(\rho_k)X_{k+1}(n),
\end{eqnarray}
where the expression of $\mathcal{Q}_i$ depends on the expression of the observable (diagonal or not).

In the non diagonal case in \cite{Diffusion}, an essential result is the proof of the convergence
\begin{equation}
W_n(t)=\frac{1}{\sqrt{n}}\sum_{k=1}^{[nt]}X_{k}(n)\mathop{\longrightarrow}_{n\rightarrow\infty}^{\mathcal{D}}W_t
\end{equation}
where $\mathcal{D}$ denotes the convergence in distribution and $(W_t)$ is a
standard Brownian motion. Next by using a theorem of convergence for
stochastic integrals due to Kurtz and Protter (cf \cite{WLT},\cite{OWZ}), it was shown that discrete quantum trajectories converges to solution of diffusive Belavkin equations. The convergence of $(W_n(t))$, in the non diagonal case, is proved independently of the convergence of $(\rho_{[nt]})$. Besides this convergence allows to conclude to the convergence of $(\rho_{[nt]})$ with Theorem of Kurtz and Protter (such result needs namely the convergence of the driving process before to consider the convergence of quantum trajectories).

In the diagonal case, we expect that the discrete quantum trajectory converges to the solution of jump-equation. In this case, a similar result for driving processes would be\begin{equation}N_n(t)=\sum_{k=1}^{[nt]}X_{k}(n)\mathop{\longrightarrow}_{n\rightarrow\infty}^{\mathcal{D}}\tilde{N}_t-
\int_0^tTr[\mathcal{J}(\rho_{s-})]ds.
\end{equation}
Apparently such result would need first the result of convergence of $(\rho_{[nt]})$ to the solution of jump-equation. Actually in the same spirit of definition of solution for the jump-equation involving the simultaneous existence of $(\tilde{N}_t)$ and $(\rho_t)$ because of the dependance of the two processes, in order to prove the final convergence result, we have to prove the simultaneous convergence of $(N_n(t))$ and $(\rho_{[nt]})$. Since we cannot show the convergence of $N_n(t)$, the convergence result of Kurtz and Protter cannot be applied. This justifies the use of random coupling method in order to compare in a simultaneous way all the processes.

Before to deal with this method, we present the result of the convergence of Euler scheme.

\subsection{Euler-scheme for jump-Belavkin equation}

The literature abounds in references about Euler scheme
approximation for stochastic differential equations (cf
\cite{OTS},\cite{ESD},\cite{EM}). The non-usual case of
jump-Belavkin equations is not really treated, that is why we
present the different result in this situation.

When we want to study Euler scheme
approximation for stochastic differential equations, an important property is the
Lipschitz character of the coefficients. Remind that our
equation is of the following form
\begin{eqnarray}\label{jump-equation}
\mu_t&=&\mu_0+\int_0^tf(\mu_{s-})ds\nonumber\\&&+\int_0^t\int_{[0,K]}\left[\frac{\mathcal{J}(\mu_{s-})}{Tr[\mathcal{J}(\mu_{s-})]}-\mu_{s-}\right]\mathbf{1}_{0\leq
x\leq Tr[\mathcal{J}(\mu_{s-})]}N(.,dx,ds)
\end{eqnarray}

We must transform this equation to have Lipschitz property. We want to write it in the following way
\begin{equation}
\mu_t=\mu_0+\int_0^tf(\mu_{s-})ds+\int_0^t\int_{[0,K]}\left[q(\mu_{s-})\right]\mathbf{1}_{0\leq
x\leq Tr[\mathcal{J}(\mu_{s-})]}N(.,dx,ds)
\end{equation}
where $f$ and $q$ are Lipschitz functions and defined for all
matrices. The consideration about the Lipschitz property is a pure technical
aspect and can be admit by the reader. The Euler scheme is given
by the formula $(\ref{eulerscheme})$ below.

Concerning $f$, we have seen that the solution of $(\ref{jump-equation})$ is obtained by truncation method because $f$ is not
Lipschitz but $C^\infty$. It was shown that the truncation is unnecessary since the solution is a process valued in the set of
states. As a consequence we can consider that the function is truncated and then Lipschitz. We denote by $F$ its Lipschitz
coefficient.

Concerning $q$, we must control the function defined on the states
by
\begin{equation}
g:\rho\longrightarrow\left[\frac{\mathcal{J}(\rho)}{Tr[\mathcal{J}(\rho)]}-\rho\right]\mathbf{1}_{0<
Tr[\mathcal{J}(\rho)]}
\end{equation} We transform the expression and define a function $q$ which is $C^\infty$ and such that
$$g(\rho)=q(\rho)\mathbf{1}_{0< Tr[\mathcal{J}(\rho)]}.$$
To construct the function $q$, it depends on the invertible character of $C$.

If $C$ is invertible, the function defined on the set of states by $\rho\rightarrow Tr[\mathcal{J}(\rho)]$ is continuous. With the fact
that for all state $\rho$, we have $Tr[\mathcal{J}(\rho)]>0$ and a compactness argument, the function $\rho\longrightarrow\mathcal{J}(\rho)/Tr[\mathcal{J}(\rho)]$
is extendible by a function $C^\infty$ defined for all matrices.

If $C$ is not invertible, there exists a unitary-operator $V$ and
two complex scalars $\alpha$ et $\beta$ such that
\begin{equation}\label{unit}VCV^\star=\left(\begin{array}{cc} \alpha & \beta\\
 0& 0
\end{array}\right).\end{equation}
Before to go further, we have to show that the jump-equation is equivalent under unitary modification. For this, define for any unitary-operator $V$
\begin{eqnarray}
\mathcal{J}_V(\rho)&=&VCV^\star(\rho)(VCV^\star)^\star\\
f_V(\rho)&=&-i[VHV^\star,\rho]-\frac{1}{2}\left\lbrace
VCV^\star(VCV^\star)^\star,\rho\right\rbrace\nonumber\\&&+VCV^\star\rho
(VCV^\star)^\star-\mathcal{J}_V(\rho)+Tr[\mathcal{J}_V(\rho)]\rho\\
g_V(\rho)&=&\left[\frac{\mathcal{J_V}(\rho)}{Tr[\mathcal{J_V}(\rho)]}-\rho\right]\mathbf{1}_{0<
Tr[\mathcal{J}_V(\rho)]}.
\end{eqnarray}
We have the following proposition which expresses unitary equivalence.
\begin{pr}\label{Unitarytransfo}
Let $V$ be any unitary operator and let $(\mu_t)$ be the solution
of the jump Belavkin equation, then the process $(\gamma_t:=V\mu_t
V^{\star})$ valued on the set of states satisfies\begin{equation}
\gamma_t=\gamma_0+\int_0^tf_V(\gamma_{s-})ds+\int_0^t\int_{[0,K]}\left[g_V(\gamma_{s-})\right]\mathbf{1}_{0\leq
x\leq Tr[\mathcal{J}_V(\gamma_{s-})]}N(.,dx,ds).
\end{equation}
\end{pr}

The proof is a straightforward computation. Such
unitary equivalence allows us to transform $g$ without change Lipschitz
property of $f$.

Now let see how we can construct the function $q$ in the case where $C$ is not invertible. Let $V$ be the unitary operator involved in expression $(\ref{unit})$, we get
$$g_V(\rho)=\left[\left(\begin{array}{cc} 1 & 0\\
 0& 0
\end{array}\right)-\rho\right]\mathbf{1}_{0< Tr[\mathcal{J}_V(\rho)]}.$$
Hence, the expression of $q$ is clear and by using the unitary
transformation given by Proposition 4, we can consider the following equation\begin{equation}
\mu_t=\mu_0+\int_0^tf(\mu_{s-})ds+\int_0^t\int_{[0,1]}[q(\mu_{s-})]\mathbf{1}_{0\leq
x\leq Tr[\mathcal{J}(\mu_{s-})]}N(.,dx,ds)
\end{equation}
which admits a unique solution by Theorem $4$.

Another time, the fact that $q$ is $C^\infty$ allows to consider that it is a Lipschitz function with a truncature method. We
denote by $Q$ the Lipschitz coefficient of $q$.
\bigskip

With these technical precautions, we can consider the Euler scheme
\begin{equation}\label{eulerscheme}\theta_{k+1}=\theta_k+\frac{1}{n}f(\theta_k)+
\int_{\frac{k}{n}}^{\frac{k+1}{n}}\int_{[0,K]}[q(\theta_k)]\mathbf{1}_{0\leq
x\leq Re(Tr[\mathcal{J}(\theta_{k})])}N(.,dx,ds).
\end{equation}
Let fixe an interval $[0,T]$ and for all $t<T$, we define
$k_t=\max\{k\in\{0,1,\ldots\}/\frac{k}{n}\leq t\}$. For all
$t$ in $]\frac{k}{n},\frac{k+1}{n}]$, we put
\begin{equation}\label{EULERSCHEME}
\tilde{\theta}_t=\theta_k+\int_{\frac{k}{n}}^tf(\theta_k)ds+\int_{\frac{k}{n}}^t\int_{[0,1]}[q(\theta_k)]\mathbf{1}_{0\leq
x\leq Re(Tr[\mathcal{J}(\theta_{k})])}N(.,dx,ds).
\end{equation}
It is worth noticing that we have $\tilde{\theta}_{\frac{k}{n}}=\theta_k$ for all $k$. We
then have for $t<T$
\begin{eqnarray*}
\tilde{\theta}_t(n)&=&\mathop{\sum}_{k=0}^{k_t-1}\int_{\frac{k}{n}}^{\frac{k+1}{n}}f(\theta_k)ds+\mathop{\sum}_{k=0}^{k_t-1}
\int_{\frac{k}{n}}^{\frac{k+1}{n}}\int_{[0,K]}[q(\theta_k)]\mathbf{1}_{0\leq
x\leq
Re(Tr[\mathcal{J}(\theta_{k})])}N(.,dx,ds)\nonumber\\
&&+\int_{k_t}^{t}f(\theta_{k_t})ds+
\int_{k_t}^{t}\int_{[0,K]}[q(\theta_{k_t})]\mathbf{1}_{0\leq x\leq
Re(Tr[\mathcal{J}(\theta_{k_t})])}N(.,dx,ds).
\end{eqnarray*}
Likewise for the solution $(\mu_t)$ of the Belavkin equation, we can write
\begin{eqnarray}
\mu_t&=&\mathop{\sum}_{k=0}^{k_t-1}\int_{\frac{k}{n}}^{\frac{k+1}{n}}f(\mu_{s-})ds+\mathop{\sum}_{k=0}^{k_t-1}
\int_{\frac{k}{n}}^{\frac{k+1}{n}}\int_{[0,1]}[q(\mu_{s-})]\mathbf{1}_{0\leq
x\leq
Tr[\mathcal{J}(\mu_{s-})]}N(.,dx,ds)\nonumber\\
&&+\int_{k_t}^{t}f(\mu_{s-})ds+
\int_{k_t}^{t}\int_{[0,K]}[q(\mu_{s-})]\mathbf{1}_{0\leq x\leq
Tr[\mathcal{J}(\mu_{s-})]}N(.,dx,ds).
\end{eqnarray}

Before to express the convergence theorem, we need the following
proposition.
\begin{pr}\label{continuity}
Let $(\mu_t)$ be the solution of the jump-Belavkin equation, then
there exists a constant $M$ such that for all
$(s,t)\in\mathbf{R}_+^2$
\begin{equation}
\mathbf{E}[\Vert\mu_t-\mu_s\Vert]\leq M\vert t-s\vert.
\end{equation}
\end{pr}
\begin{pf}
From the fact that the solution of the jump-Belavkin equation is
valued on the set of states, we have for all $t>0$, the quantity
$\Vert\mu_t\Vert\leq 1$ almost surely (we do not specify the norm
because we just need the fact that the solution is bounded). For $0<s<t$, we have
\begin{equation}
\mu_t-\mu_s=\int_s^tf(\mu_{u-})du+\int_{s}^{t}\int_{[0,1]}[q(\mu_{u-})]\mathbf{1}_{0\leq
x\leq Tr[\mathcal{J}(\mu_{u-})]}N(.,dx,du).
\end{equation}
By using the property of random Poisson measure, in particular the
property about the intensity measure we
have for $0<s<t$
\begin{eqnarray*}
\mathbf{E}\left[\Vert\mu_t-\mu_s\Vert\right]&\leq&\mathbf{E}\left[\left\Vert\int_s^tf(\mu_{u-})du\right\Vert\right]\\&&+\mathbf{E}\left[\left\Vert\int_{s}^{t}\int_{[0,K]}[q(\mu_{u-})]\mathbf{1}_{0\leq
x\leq
Tr[\mathcal{J}(\mu_{u-})]}N(.,dx,du)\right\Vert\right]\\&\leq&\int_s^t\mathbf{E}\left[\left\Vert
f(\mu_{u-})\right\Vert\right]
du\\&&+\mathbf{E}\left[\int_{s}^{t}\int_{[0,K]}\left\Vert[q(\mu_{u-})]\mathbf{1}_{0\leq
x\leq
Tr[\mathcal{J}(\mu_{u-})]}\right\Vert N(.,dx,du)\right]\\
&\leq&\int_s^t\mathbf{E}\left[\left\Vert f(\mu_{u-})\right\Vert\right] du+\mathbf{E}\left[\int_{s}^{t}\int_{[0,K]}\left\Vert
q(\mu_{u-})\mathbf{1}_{0\leq x\leq Tr[\mathcal{J}(\mu_{u-})]}\right\Vert dxdu\right]\\&\leq&\int_s^t\left(\sup_{\Vert R\Vert\leq
1}\Vert f(R)\Vert+2\right)du\\&\leq& M\left(t-s\right)
\end{eqnarray*}
where $M$ is a constant. The result is then proved.
\end{pf}\\

We can now express the theorem concerning the convergence of the
Euler scheme. The particularity of the jump-equation comes from the stochastic intensity depending on the solution.

\begin{thm}\label{premier} Let $T>0$, let $(\tilde{\theta}_t)$ be the process $(\ref{EULERSCHEME})$
constructed by the Euler-scheme on $[0,T]$, let $(\mu_t)$ be the unique solution of the jump-Belavkin equation $(\ref{integral})$.

We define for $u<T$ and $n$ large enough
\begin{equation}Z_u(n)=\mathbf{E}\left[\mathop{\sup}_{0\leq t\leq
u}\left\Vert\tilde{\theta}_t(n)-\mu_t\right\Vert\right].
\end{equation}
So there exists a constant $\Gamma$ which is independent of $n$ such that for all $u<T$
\begin{equation}Z_u(n)\leq \frac{\Gamma}{n}.
\end{equation}
Let $\mathcal{D}\left(\left[0,T\right]\right)$ denotes the space
of c\`adl\`ag matrices processes endowed with the Skorohod topology.

Hence the Euler scheme approximation $(\tilde{\theta}_t)$
converges in distribution in
$\mathcal{D}\left(\left[0,T\right]\right)$ for all $T$ to the
process-solution $(\mu_t)$ of the jump-Belavkin equation.
\end{thm}

Before to give the proof, it is interesting to compare this result with the classic ones. In the literature, similar result appears often by using a $L_2$ norm \cite{OTS}, that is terms like $$\mathbf{E}\left[\mathop{\sup}_{0\leq t\leq
u}\left\Vert\tilde{\theta}_t(n)-\mu_t\right\Vert^2\right]$$ are usually considered (to apply It\^o Isometry result). Results like $Z_u(n)\leq \Gamma/n^2$ are then shown and almost surely convergence are then obtained. In the following proof, it is explained why such norm is not appropriate in our case.\\
\begin{pf}
The equations concerning the Euler scheme and the solution of the jump-Belavkin equation give us the following formula
\begin{eqnarray*}
\tilde{\theta}_t(n)-\mu_t&=&\mathop{\sum}_{k=0}^{k_t-1}\int_{\frac{k}{n}}^{\frac{k+1}{n}}[f(\theta_k)-f(\mu_{s-})]ds
+\int_{k_t}^{t}[f(\theta_{k_t})-f(\mu_{s-})]ds\\&&+\mathop{\sum}_{k=0}^{k_t-1}
\int_{\frac{k}{n}}^{\frac{k+1}{n}}\int_{[0,K]}\Big{(}q(\theta_k)\mathbf{1}_{0\leq
x\leq
Re(Tr[\mathcal{J}(\theta_{k})])}\\&&\hspace{3,5cm}-q(\mu_{s-})\mathbf{1}_{0\leq
x\leq Tr[\mathcal{J}(\mu_{s-})]}\Big{)}N(.,dx,ds)\\&&+
\int_{k_t}^{t}\int_{[0,K]}\Big{(}q(\theta_{k_t})\mathbf{1}_{0\leq x\leq
Re(Tr[\mathcal{J}(\theta_{k_t})])}\\&&\hspace{2,5cm}-q(\mu_{s-})\mathbf{1}_{0\leq
x\leq Tr[\mathcal{J}(\mu_{s-})]}\Big{)}N(.,dx,ds).
\end{eqnarray*}
Let consider for $u<T$, the quantity $Z_u(n)=\mathbf{E}\left[\mathop{\sup}_{0\leq t\leq
u}\left\Vert\tilde{\theta}_t(n)-\mu_t\right\Vert\right]$. Let consider separately the drift term and the term concerning the random
measure. For the drift term, by the fact that $f$ is Lipschitz, we have
\begin{eqnarray*}
&&\mathbf{E}\left[\mathop{\sup}_{0\leq t\leq
u}\mathop{\sum}_{k=0}^{k_t-1}\int_{\frac{k}{n}}^{\frac{k+1}{n}}\Vert
f(\theta_k)-f(\mu_{s-})\Vert ds+\int_{k_t}^{t}\Vert
f(\theta_{k_t})-f(\mu_{s-})\Vert
ds\right]\\
&\leq&
\mathbf{E}\left[\mathop{\sum}_{k=0}^{k_u-1}\int_{\frac{k}{n}}^{\frac{k+1}{n}}\left\Vert
f(\tilde{\theta}_{\frac{k}{n}})-f(\mu_{s-})\right\Vert
ds+\int_{k_u}^{u}\left\Vert
f(\tilde{\theta}_{\frac{k_u}{n}})-f(\mu_{s-})\right\Vert ds\right]\\
&\leq&
\mathop{\sum}_{k=0}^{k_u-1}\int_{\frac{k}{n}}^{\frac{k+1}{n}}\mathbf{E}\left[\left\Vert
f(\tilde{\theta}_{\frac{k}{n}})-f(\mu_{\frac{k}{n}})\right\Vert\right]
ds+\int_{k_u}^{u}\mathbf{E}\left[\left\Vert
f(\tilde{\theta}_{\frac{k_u}{n}})-f(\mu_{\frac{k_u}{n}})\right\Vert\right] ds\\
&&
+\mathop{\sum}_{k=0}^{k_u-1}\int_{\frac{k}{n}}^{\frac{k+1}{n}}\mathbf{E}\left[\left\Vert
f(\mu_{s-})-f(\mu_{\frac{k}{n}})\right\Vert\right]
ds+\int_{k_u}^{u}\mathbf{E}\left[\left\Vert
f(\mu_{s-})-f(\mu_{\frac{k_u}{n}})\right\Vert\right] ds\\&\leq&
\mathop{\sum}_{k=0}^{k_u-1}\int_{\frac{k}{n}}^{\frac{k+1}{n}}F\mathbf{E}\left[\left\Vert
\tilde{\theta}_{\frac{k}{n}}-\mu_{\frac{k}{n}}\right\Vert\right]
ds+\int_{k_u}^{u}F\mathbf{E}\left[\left\Vert
\tilde{\theta}_{\frac{k_u}{n}}-\mu_{\frac{k_u}{n}}\right\Vert\right] ds\\
&&
+\mathop{\sum}_{k=0}^{k_u-1}\int_{\frac{k}{n}}^{\frac{k+1}{n}}\mathbf{E}\left[\left\Vert
f(\mu_{s-})-f(\mu_{\frac{k}{n}})\right\Vert\right]
ds+\int_{k_u}^{u}\mathbf{E}\left[\left\Vert
f(\mu_{s-})-f(\mu_{\frac{k_u}{n}})\right\Vert\right] ds\\
&\leq&
\mathop{\sum}_{k=0}^{k_u-1}\int_{\frac{k}{n}}^{\frac{k+1}{n}}F\mathbf{E}\left[\sup_{0\leq
t\leq s}\left\Vert \tilde{\theta}_{t}-\mu_{t}\right\Vert\right]
ds+\int_{k_u}^{u}F\mathbf{E}\left[\sup_{0\leq t\leq s}\left\Vert
\tilde{\theta}_{t}-\mu_{t}\right\Vert\right] ds\\
&&+
\mathop{\sum}_{k=0}^{k_u-1}\int_{\frac{k}{n}}^{\frac{k+1}{n}}FM\left(s-\frac{k}{n}\right)ds+
\int_{k_u}^{u}FM\left(s-\frac{k_u}{n}\right)
ds\\
&\leq& A\left(\int_0^uZ_sds+\frac{1}{n}\right)\,\,\textrm{(A is a suitable
constant).}
\end{eqnarray*}

The analysis of the random measure terms is more complicated. Let fixe an index $k$, thanks to the properties of random measure, we
have
\begin{eqnarray}\label{eq}
&&\mathbf{E}\left[\left\Vert\int_{\frac{k}{n}}^{\frac{k+1}{n}}\int_{[0,K]}[q(\theta_k)]\mathbf{1}_{0\leq
x\leq
Re(Tr[\mathcal{J}(\theta_{k})])}-[q(\mu_{s-})]\mathbf{1}_{0\leq
x\leq
Tr[\mathcal{J}(\mu_{s-})]}N(.,dx,ds)\right\Vert\right]\nonumber\\
&\leq&
\mathbf{E}\left[\int_{\frac{k}{n}}^{\frac{k+1}{n}}\int_{[0,K]}\left\Vert[q(\theta_k)]\mathbf{1}_{0\leq
x\leq
Re(Tr[\mathcal{J}(\theta_{k})])}-[q(\mu_{s-})]\mathbf{1}_{0\leq
x\leq
Tr[\mathcal{J}(\mu_{s-})]}\right\Vert N(.,dx,ds)\right]\nonumber\\
&\leq&
\mathbf{E}\left[\int_{\frac{k}{n}}^{\frac{k+1}{n}}\int_{[0,K]}\left\Vert[q(\theta_k)]\mathbf{1}_{0\leq
x\leq
Re(Tr[\mathcal{J}(\theta_{k})])}-[q(\mu_{s-})]\mathbf{1}_{0\leq
x\leq
Tr[\mathcal{J}(\theta_{k})]}\right\Vert N(.,dx,ds)\right]\nonumber\\
&&+\mathbf{E}\left[\int_{\frac{k}{n}}^{\frac{k+1}{n}}\int_{[0,K]}\left\Vert[q(\mu_{s-})]\right\Vert\times
\left\vert\mathbf{1}_{0\leq
x\leq Tr[\mathcal{J}(\mu_{s-})]}-\mathbf{1}_{0\leq x\leq
Re(Tr[\mathcal{J}(\theta_{k})])}\right\vert N(.,dx,ds)\right]\nonumber\\
&\leq&
\mathbf{E}\left[\int_{\frac{k}{n}}^{\frac{k+1}{n}}\int_{[0,K]}\left\Vert
q(\tilde{\theta}_{\frac{k}{n}})-q(\mu_{s-})\right\Vert
N(.,dx,ds)\right]\nonumber\\
&&+2\mathbf{E}\left[\int_{\frac{k}{n}}^{\frac{k+1}{n}}\int_{[0,K]}\left\vert\mathbf{1}_{0\leq
x\leq Tr[\mathcal{J}(\mu_{s-})]}-\mathbf{1}_{0\leq x\leq
Tr[\mathcal{J}(\theta_{k})]}\right\vert
N(.,dx,ds)\right].\end{eqnarray}  As $q$ is bounded by $2$ on the
set of states, we have
\begin{eqnarray}
(\ref{eq})&\leq&
\mathbf{E}\left[\int_{\frac{k}{n}}^{\frac{k+1}{n}}\int_{[0,K]}Q\left\Vert\tilde{\theta}_{\frac{k}{n}}-\mu_{s-}\right\Vert
N(.,dx,ds)\right]\nonumber\\
&&+2\mathbf{E}\Bigg{[}\int_{\frac{k}{n}}^{\frac{k+1}{n}}\int_{[0,K]}\Big{(}\mathbf{1}_{0\leq
x\leq \max(Tr[\mathcal{J}(\mu_{s-})],
Re(Tr[\mathcal{J}(\tilde{\theta}_{\frac{k}{n}})]))}\nonumber\\&&\hspace{3,5cm}-\mathbf{1}_{0\leq
x\leq \min(Tr[\mathcal{J}(\mu_{s-})],
Re(Tr[\mathcal{J}(\tilde{\theta}_{\frac{k}{n}})]))}\Big{)}
N(.,dx,ds)\Bigg{]}.\nonumber\end{eqnarray}
Furthermore, we have
$Re(Tr[\mathcal{J}(\mu_{s-})])=Tr[\mathcal{J}(\mu_{s-})]$ for all $s$. Hence, by linearity and continuity of the function trace, for
any matrices $A$ and $B$, there exists a constant $R$ such that $$\left\vert Re(
Tr[\mathcal{J}(A)])-Re(Tr[\mathcal{J}(B)])\right\vert\leq R\left\Vert A-B\right\Vert.$$ It implies
\begin{eqnarray}
(\ref{eq})&\leq&
\mathbf{E}\left[\int_{\frac{k}{n}}^{\frac{k+1}{n}}Q\left\Vert\tilde{\theta}_{\frac{k}{n}}-\mu_{s-}\right\Vert
ds\right]+\mathbf{E}\left[\int_{\frac{k}{n}}^{\frac{k+1}{n}}\left\vert
Re(
Tr[\mathcal{J}(\tilde{\theta}_{\frac{k}{n}})])-Tr[\mathcal{J}(\mu_{s-})]\right\vert
ds\right]\nonumber\\&\leq&
\int_{\frac{k}{n}}^{\frac{k+1}{n}}(R+Q)\mathbf{E}\left[\left\Vert\tilde{\theta}_{\frac{k}{n}}-\mu_{s-}\right\Vert\right]
ds\nonumber\\&\leq&
\int_{\frac{k}{n}}^{\frac{k+1}{n}}\left(R+Q\right)\mathbf{E}\left[\left\Vert\tilde{\theta}_{\frac{k}{n}}
-\mu_{\frac{k}{n}}\right\Vert\right]ds+
\int_{\frac{k}{n}}^{\frac{k+1}{n}}\left(R+Q\right)\mathbf{E}\left[\left\Vert\tilde{\mu}_{\frac{k}{n}}-\mu_{s-}\right\Vert\right]
ds\nonumber\\
&\leq&
(R+Q)\left(\int_{\frac{k}{n}}^{\frac{k+1}{n}}Z_sds+\frac{1}{n^2}\right)\,\,\,\textrm{(we
do the same as the drift term).}
\end{eqnarray}

The term between $k_{u}$ and $u$ can be treated in the same way.
By summing, we obtain finally the same type of inequality for the
term with the random measure. As a consequence there exist two
constants $F_1$ and $F_2$ which depend only on $T$ such that
\begin{equation}
Z_u\leq F_1\int_0^uZ_sds+\frac{F_2}{n}.
\end{equation}
The Gronwall Lemma implies that there exists a constant $\Gamma$ such that for all $u<T$
\begin{equation}
Z_u(n)\leq \frac{\Gamma}{n}
\end{equation}
where $\Gamma$ is a constant independent of $n$. The convergence
in $\mathcal{D}\left([0,T]\right)$ is an easy consequence of the
above inequality. The result is then proved.

Now we can justify why the $L_2$ norm is not appropriate to deal with such equation. Indeed, the last term of (\ref{eq}) involves an integral with a difference of two indicator functions. The difference of the indicator functions is equal to zero or one and it gives the same result if we choose the $L_2$ norm, that is
$$\left\vert\mathbf{1}_{0\leq
x\leq Tr[\mathcal{J}(\mu_{s-})]}-\mathbf{1}_{0\leq x\leq
Tr[\mathcal{J}(\theta_{k})]}\right\vert^2=\left\vert\mathbf{1}_{0\leq
x\leq Tr[\mathcal{J}(\mu_{s-})]}-\mathbf{1}_{0\leq x\leq
Tr[\mathcal{J}(\theta_{k})]}\right\vert.
$$
After in the above proof, the integral of this term is calculated and if we use $L_2$ norm, we loose the homogeneity in term of $L_2$ norm.  As the final result
relies on Gronwall lemma lemma, this homogeneity is actually necessary to obtain an appropriate estimation.
\end{pf}\\

In the following section, we compare the discrete process with the Euler scheme.

\subsection{Convergence of the discrete process}

This section is devoted to the random coupling method of the discrete quantum trajectory and the continuous quantum trajectory.

Consider the probability space $(\Omega,\mathcal{F},\mathcal{F}_t,P)$ where the solution $(\mu_t)$ of the jump Belavkin equation
$(\ref{integral})$ is defined. Let construct the discrete quantum trajectory in this space. Let $n$ be fixed, we define the following
sequence of random variable which are defined on the set of states
 \begin{equation}\tilde{\nu}_{k+1}(\eta,\omega)=\mathbf{1}_{N(\omega,G_k(\eta))>0}
 \end{equation}
 where $G_k(\eta)=\left\{(t,u)/\frac{k}{n}\leq t<\frac{k+1}{n},0\leq
u\leq-n\ln(Tr[\mathcal{L}_{0}(n)(\eta)])\right\}$.

 Let $\rho_0=\rho$ be any
state, we define the process $(\tilde{\rho}_k)$ for $k<[nT]$ by the recursive
formula
\begin{eqnarray}\label{eqq}
\tilde{\rho}_{k+1}&=&\mathcal{L}_{0}(\tilde{\rho}_k)+\mathcal{L}_{1}(\tilde{\rho}_k)\nonumber\\
&&+\left[-\frac{\mathcal{L}_{0}(\tilde{\rho}_k)}{Tr[\mathcal{L}_{0}(\tilde{\rho}_k)]}+
\frac{\mathcal{L}_{1}(\tilde{\rho}_k)}{Tr[\mathcal{L}_{1}(\tilde{\rho}_k)]}\right]
\left(\tilde{\nu}_{k+1}(\tilde{\rho}_k,.)-Tr[\mathcal{L}_{1}(\tilde{\rho}_k)]\right).
\end{eqnarray}
This random sequence and the operators $\mathcal{L}_{i}(n)$
depend naturally
 on $n$ following the asymptotic of the unitary evolution. Thanks to the Poisson distribution property, the following
proposition is obvious.

\begin{pr}Let $T>0$ be fixed. The discrete process $(\tilde{\rho}_k)_{k<[nT]}$ defined by $(\ref{eqq})$
has the same distribution of the discrete quantum trajectory $(\rho_k)_{k<[nT]}$ defined by the quantum repeated measurements
principle.
\end{pr}

 This proposition is a consequence of the fact that for all Borel subset $B\in\mathcal{B}(\mathbf{R}^2)$, we get
 $$P[N(B)=k]=\frac{\Lambda(B)^k}{k!}\exp(-\Lambda(B)),$$
 where $\Lambda$ denotes the Lebesgue measure.

 In $(\Omega,\mathcal{F},\mathcal{F}_t,P)$,  the process
$(\tilde{\rho}_k)$ satisfies the same asymptotic than the discrete quantum trajectory, that is,
 \begin{equation}\label{discrete1}
\tilde{\rho}_{k+1}-\tilde{\rho}_k=\frac{1}{n}[f(\tilde{\rho}_k)+\circ_{\tilde{\rho}_k}(1)]+\left[\frac{\mathcal{J}(\tilde{\rho}_k)}{Tr(\mathcal{J}(\tilde{\rho}_k))}-\tilde{\rho}_k+\circ_{\tilde{\rho}_k}(1)\right]
\tilde{\nu}_{k+1}(\tilde{\rho}_k,.).
\end{equation}

Before to compare the discrete process $(\ref{eqq})$, we need
another process which concerns the approximation of the intensity. In $(\Omega,\mathcal{F},\mathcal{F}_t,P)$, we
define the random variable sequence defined on the set of states
\begin{equation}\overline{\nu}_{k+1}(\eta,\omega)=\mathbf{1}_{N(\omega,H_k(\eta))>0}
 \end{equation}
 where $H_k(\eta)=\{(t,u)/\frac{k}{n}\leq t<\frac{k+1}{n},0\leq u\leq
Tr[\mathcal{J}(\eta)]\}$.
  Let
 $\overline{\rho}_0=\rho$ be a state, we define the following process in
 $(\Omega,\mathcal{F},\mathcal{F}_t,P)$, for $k<[nT]$, we put
 \begin{eqnarray}\label{eqqq}
\overline{\rho}_{k+1}&=&\mathcal{L}_{0}(\overline{\rho}_k)+\mathcal{L}_{1}(\overline{\rho}_k)\nonumber\\
&&+\left[-\frac{\mathcal{L}_{0}(\overline{\rho}_k)}{Tr[\mathcal{L}_{0}(\overline{\rho}_k)]}+
\frac{\mathcal{L}_{1}(\overline{\rho}_k)}{Tr[\mathcal{L}_{1}(\overline
{\rho}_k)]}\right]
\left(\overline{\nu}_{k+1}(\overline{\rho}_k,.)-Tr[\mathcal{L}_{1}(\overline{\rho}_k)]\right).
\end{eqnarray}
Hence, we have the same asymptotic form
 \begin{equation}\label{form}\overline{\rho}_{k+1}-\overline{\rho}_k=\frac{1}{n}[f(\overline{\rho}_k)+\circ_{\overline{\rho}_k}(1)]+\left[\frac{\mathcal{J}(\overline{\rho}_k)}{Tr(\mathcal{J}(\overline{\rho}_k))}-\overline{\rho}_k+\circ_{\overline{\rho}_k}(1)\right]
\overline{\nu}_{k+1}(\overline{\rho}_k,.).
\end{equation}
Regarding the process $(\tilde{\rho}_k)_{0\leq k\leq [nT]}$ and $(\overline{\rho}_k)_{0\leq k\leq [nT]}$, we have the following proposition.

\begin{pr}\label{P2} Let $(\tilde{\rho}_k)_{0\leq k\leq [nT]}$ be the discrete quantum trajectory
defined by the formula $(\ref{discrete1}$ and let $(\overline{\rho}_k)_{0\leq k\leq [nT]}$ be the sequence defined
by the formula $(\ref{eqqq})$. Let assume that the two sequences are defined by the same initial state $\rho$.

 For $k\leq[nT]$, we define $$A_k(n)=\mathbf{E}\left[\sup_{0<i\leq
k}\left\Vert\tilde{\rho}_i(n)-\overline{\rho}_i(n)\right\Vert\right].$$
We have for all $k\leq[nT]$
$$A_k(n)\leq\circ\left(\frac{1}{n}\right)$$
where the little $\circ$ is uniform in $k$.
\end{pr}
\begin{pf}
Remind that the discrete quantum trajectory satisfies
\begin{equation}
\tilde{\rho}_{k+1}-\tilde{\rho}_k=\frac{1}{n}[f(\tilde{\rho}_k)+\circ_{\tilde{\rho}_k}(1)]+
\left[\frac{\mathcal{J}(\tilde{\rho}_k)}{Tr[\mathcal{J}(\tilde{\rho}_k)]}-\tilde{\rho}_k+\circ_{\tilde{\rho}_k}(1)\right]
\tilde{\nu}_{k+1}(\tilde{\rho}_k,.).
\end{equation}
We can remark that all the rest $\circ_{\tilde{\rho}^k}(1)$ are
uniform in $k$ because the process $(\rho_k)$ is valued on the set
of states and so is bounded. Hence, we can write this equation in
the following way using $f$ and $q$
\begin{equation}
\tilde{\rho}_{k+1}-\tilde{\rho}_k=\frac{1}{n}[f(\tilde{\rho}_k)+
\circ_{\tilde{\rho}_k}(1)]+[q(\tilde{\rho}_k)+\circ_{\tilde{\rho}_k}(1)]\tilde{\nu}_{k+1}(\tilde{\rho}_k,.).
\end{equation}
We have a similar asymptotic form for the process $(\ref{eqqq})$. As a consequence, we can compare the two processes
\begin{eqnarray*}
\tilde{\rho}_i-\overline{\rho}_i&=&\sum_{j=0}^{i-1}\left[\frac{1}{n}(f(\tilde{\rho}_j)-f(\overline{\rho}_j)
+\circ_{\tilde{\rho}_k}(1)-\circ_{\overline{\rho}_k}(1)\right]\\&&+\sum_{j=0}^{i-1}\Big[\left(q(\tilde{\rho}_j)
+\circ_{\tilde{\rho}_j}(1)\right)\tilde{\nu}_{j+1}(\tilde{\rho}_j,.)-\left(q(\overline{\rho}_j)
+\circ_{\overline{\rho}_k}(1)\right)\overline{\nu}_{j+1}(\overline{\rho}_j,.)\Big].
\end{eqnarray*}
Hence, we have
\begin{eqnarray}\label{supp}
\sup_{0<i\leq
k}\left\Vert\tilde{\rho}_i-\overline{\rho}_i\right\Vert&\leq&\sum_{j=0}^{k-1}\frac{1}{n}\Big\Vert(f(\tilde{\rho}_j)
-f(\overline{\rho}_j)+\circ_{\tilde{\rho}^k}(1)\Big\Vert\nonumber\\&&+\sum_{j=0}^{k-1}\Big\Vert
(q(\tilde{\rho}_j)+\circ_{\tilde{\rho}_j}(1))\tilde{\nu}_{j+1}(\tilde{\rho}_j,.)
-q(\overline{\rho}_j)\overline{\nu}_{j+1}(\overline{\rho}_j,.)\Big\Vert\nonumber\\
&\leq&\sum_{j=0}^{k-1}\frac{1}{n}F\Big\Vert\tilde{\rho}_j-\overline{\rho}_j\Big\Vert
+\sum_{j=0}^{k-1}\Big\Vert
(q(\tilde{\rho}_j)+\circ_{\tilde{\rho}_j}(1))(\tilde{\nu}_{j+1}(\tilde{\rho}_j,.)-\overline{\nu}_{j+1}(\overline{\rho}_j,.)\Big\Vert
\nonumber\\&&+\sum_{j=0}^{k-1}\Big\Vert
q(\tilde{\rho}_j)+\circ_{\tilde{\rho}_j}(1))-q(\overline{\rho}_j)-\circ_{\overline{\rho}_j}(1)))\overline{\nu}_{j+1}(\overline{\rho}_j,.)\Big\Vert
\end{eqnarray}
By defining the filtration
$\mathcal{G}_j=\sigma\{\tilde{\nu}_{k}(\tilde{\rho}_{k-1},.),\overline{\nu}_{k}(\overline{\rho}_{k-1},.),0<k\leq
j\}$ for $j>0$, we have by the independence of the increments of a
Poisson process
\begin{eqnarray}\label{supa}
&&\mathbf{E}\left[\Vert
q(\tilde{\rho}^j)+\circ_{\tilde{\rho}_j}(1))-q(\overline{\rho}_j)-\circ_{\overline{\rho}_j}(1)))\overline{\nu}_{j+1}(\overline{\rho}_j,.)\Vert\right]\nonumber\\
&=&\mathbf{E}\left[\mathbf{E}\left[\Vert
q(\tilde{\rho}_j)+\circ_{\tilde{\rho}_j}(1))-q(\overline{\rho}_j)-\circ_{\overline{\rho}_j}(1)))\overline{\nu}_{j+1}(\overline{\rho}_j,.)\Vert/\mathcal{G}_j\right]\right]\nonumber\\
&=&\mathbf{E}\left[\Vert
q(\tilde{\rho}_j)+\circ_{\tilde{\rho}_j}(1))-q(\overline{\rho}_j)-\circ_{\overline{\rho}_j}(1)))\mathbf{E}\left[\overline{\nu}_{j+1}(\overline{\rho}_j,.)\Vert/\mathcal{G}_j\right]\right]\nonumber\\
&=&\mathbf{E}\left[\Vert
q(\tilde{\rho}_j)+\circ_{\tilde{\rho}_j}(1)-q(\overline{\rho}_j)-\circ_{\overline{\rho}_j}(1)\Vert\right]\left(1-exp\left(-\frac{1}{n}Tr[\mathcal{J}(\overline{\rho}_j)]\right)\right)\nonumber\\
&\leq&\mathbf{E}\left[Q\Vert
\tilde{\rho}_j-\overline{\rho}_j\Vert\right]\left(\frac{1}{n}Tr[\mathcal{J}(\overline{\rho}_j)]+\circ\left(\frac{1}{n}\right)\right)+\circ\left(\frac{1}{n}\right),
\end{eqnarray}
because all the rest are uniform in $j$. The same way, by using the filtration, we get
\begin{eqnarray*}
&&\mathbf{E}\Big[\left\Vert
(q(\tilde{\rho}_j)+\circ_{\tilde{\rho}_j}(1))(\tilde{\nu}_{j+1}(\tilde{\rho}_j,.)-\overline{\nu}_{j+1}(\overline{\rho}_j,.)\right\Vert\Big]\\
&=&\mathbf{E}\Big[\left\Vert
q(\tilde{\rho}_j)+\circ_{\tilde{\rho}_j}(1))\right\Vert\mathbf{E}\left[\left\vert\tilde{\nu}_{j+1}(\tilde{\rho}_j,.)-\overline{\nu}_{j+1}(\overline{\rho}_j,.)\right\vert/\mathcal{G}_j\right]\Big].
\end{eqnarray*}
For the second term of the product, by definition, we have
\begin{eqnarray*}
&&\mathbf{E}\Big[\left\vert\tilde{\nu}_{j+1}(\tilde{\rho}_j,.)-\overline{\nu}_{j+1}(\overline{\rho}_j,.)\right\vert/\mathcal{G}_j\Big]\\
&=&\mathbf{E}\Big[\left\vert\mathbf{1}_{N(.,G_{j}(\tilde{\rho}_j))>0}-\mathbf{1}_{N(.,H_{j}(\overline{\rho}_j))>0}\right\vert/\mathcal{G}_j\Big]\\
&=&\mathbf{E}\left[\mathbf{1}_{\{N(.,G_{j}(\tilde{\rho}_j))>0\}\bigtriangleup\{N(.,H_{j}(\overline{\rho}_j))>0\}}/\mathcal{G}_j\right]\\
&=&P\left[\{N(.,G_{j}(\tilde{\rho}_j))>0\}\bigtriangleup\{N(.,H_{j}(\overline{\rho}_j))>0\}/\mathcal{G}_j\right].\\
\end{eqnarray*}
We denote by $W_j=\{(t,u)/\frac{j}{n}\leq
t<\frac{j+1}{n},\min\left(Tr[\mathcal{J}(\overline{\rho}_j)],-n\ln(Tr[\mathcal{L}_0(\rho_j)])\right)\leq
u\leq\max\left(Tr[\mathcal{J}(\overline{\rho}_j)],-n\ln(Tr[\mathcal{L}_0(\rho_j)])\right)\}.$
Hence we have
\begin{eqnarray*}
\mathbf{E}\Big[\left\vert\tilde{\nu}_{j+1}(\tilde{\rho}_j,.)-\overline{\nu}_{j+1}(\overline{\rho}_j,.)\right\vert/\mathcal{G}_j\Big]
&=&P[\mathbf{1}_{W_j>0}/\mathcal{G}_j]\\&=&1-\exp\bigg(-\frac{1}{n}\Big(\min\left(Tr[\mathcal{J}(\overline{\rho}_j)],-n\ln(Tr[\mathcal{L}_0(\rho_j)])\right)\\&&\hspace{3
cm}-\max\left(\ldots,-n\ln(Tr[\mathcal{L}_0(\rho_j)])\right)\Big)\bigg)\\
&=&1-\exp\left(-\frac{1}{n}\Big\vert
Tr[\mathcal{J}(\overline{\rho}^j)]+n\ln[Tr[\mathcal{L}_0(\tilde{\rho}_j)]\Big\vert\right)\\
&=&\frac{1}{n}\Big\vert
Tr[\mathcal{J}(\overline{\rho}_j)]+n\ln\left(Tr[\mathcal{L}_0(\tilde{\rho}_j)]\right)\Big\vert+\circ\left(\frac{1}{n}\right).
\end{eqnarray*}
Besides we have
$Tr[\mathcal{L}_0(\tilde{\rho}_j)]=p_{j+1}=1-\frac{1}{n}Tr[\mathcal{J}(\rho_j)]+\circ(\frac{1}{n})$,
hence
\begin{eqnarray}\label{supb}
\mathbf{E}\left[\Big\vert\tilde{\nu}_{j+1}(\tilde{\rho}_j,.)-\overline{\nu}_{j+1}(\overline{\rho}_j,.)\Big\vert/\mathcal{G}_j\right]&=&\frac{1}{n}\Big\vert
Tr[\mathcal{J}(\overline{\rho}_j)]-Tr[\mathcal{J}(\tilde{\rho}_j)]\Big\vert+\circ\left(\frac{1}{n}\right).
\end{eqnarray}
As $(\tilde{\rho}_k)$ is a process valued in the set of state, it is uniformly bounded, we then have
\begin{equation}
 \mathbf{E}\Big[\left\Vert
(q(\tilde{\rho}_j)+\circ_{\tilde{\rho}_j}(1))(\tilde{\nu}_{j+1}(\tilde{\rho}_j,.)-\overline{\nu}_{j+1}(\overline{\rho}_j,.)\right\Vert\Big]\leq K\mathbf{E}\Big[\Vert\overline{\rho}_j-\tilde{\rho}_j\Vert\Big]+\circ\left(\frac{1}{n}\right).
\end{equation}
By taking expectation in $(\ref{supp})$ and using the inequalities $(\ref{supa},\ref{supb})$, we obtain finally the inequality
\begin{eqnarray}A_k&\leq&
\sum_{j=0}^{k-1}\frac{L}{n}\mathbf{E}\Big[\left\Vert\tilde{\rho}_j-\overline{\rho}_j\right\Vert\Big]+\circ\left(\frac{1}{n}\right)\nonumber\\
&\leq&\sum_{j=0}^{k-1}\frac{L}{n}A_j+\circ\left(\frac{1}{n}\right).
\end{eqnarray}
The result follows with a discrete Gronwall Lemma.
\end{pf}

Now, we can compare the process obtained by the Euler scheme and
the process defined by the formula $(\ref{form})$. The result is resumed in the
following proposition.

\begin{pr}
Let $(\overline{\rho}_k)_{0\leq k\leq [nT]}$ be the process defined by the formula $(\ref{form})$
and let $(\theta_k)_{0\leq k\leq [nT]}$ be the process obtained by the Euler scheme
of the jump-Belavkin equation. Let assume that the two sequences are defined by the same initial state $\rho$.

For $k\leq[nT]$, we define $$S_k(n)=\mathbf{E}\left[\sup_{0\leq i\leq k}\Vert\theta_i(n)-\overline{\rho}_i(n)\Vert\right].$$ Hence
there exists a constant $F$ such that for all $k\leq[nT]$
$$S_k(n)\leq\frac{F}{n}.$$
\end{pr}

The proof is based on the Gronwall Lemma but it uses
finer property of the random measure induced by the Poisson point process. This is a generalization of the Poisson approximation studied by Brown in \cite{MR704559}.\\
\begin{pf}
Thanks to the fact that random sequence $(\overline{\rho}_j)$ is bounded, the $\circ_{\overline{\rho}_j}=\circ(1)$. It implies that for
$i\leq k \leq[nT]$
\begin{eqnarray*}
\theta_i-\overline{\rho}_i&=&\sum_{j=0}^{i-1}\frac{1}{n}[f(\theta_j)+\circ(1)-f(\overline{\rho}_j)]+\sum_{j=0}^{i-1}\int_{\frac{j}{n}}^{\frac{j+1}{n}}\int_{[0,1]}[q(\theta_j)]\mathbf{1}_{0\leq
x\leq
Tr[\mathcal{J}(\theta_j)}N(.,dx,ds)\\&&-\sum_{j=0}^{i-1}[q(\overline{\rho}_j)+\circ(1)]\overline{\nu}_{j+1}(\overline{\rho}_j,.).
\end{eqnarray*}
We treat the random measure part and the drift term part
separately. Let us denote $S_k=\mathbf{E}\left[\sup_{0\leq i\leq
k}\Vert\theta_i-\overline{\rho}_i\Vert\right]$, we have
\begin{eqnarray*}
S_k&\leq&\mathbf{E}\left[\sum_{j=0}^{k-1}\frac{1}{n}\left\Vert[f(\theta_j)-f(\overline{\rho}_j)+\circ(1)]\right\Vert\right]\\
&&+\mathbf{E}\left[\sum_{j=0}^{k-1}\left\Vert\int_{\frac{j}{n}}^{\frac{j+1}{n}}\int_{[0,1]}q(\theta_j)\mathbf{1}_{0\leq x\leq
Tr[\mathcal{J}(\theta_j)]}N(.,dx,ds)\right.\right.-\left(q(\overline{\rho}_j)+\circ(1)\right)\overline{\nu}_{j+1}(\overline{\rho}_j,.)]\Bigg{\Vert}\Bigg{]}.
\end{eqnarray*}
As $f$ is Lipschitz, we have
$$\mathbf{E}\left[\sum_{j=0}^{k-1}\frac{1}{n}\left\Vert f(\theta_j)-f(\overline{\rho}_j)+\circ(1)\right\Vert\right]\leq
F\sum_{j=0}^{k-1}\frac{1}{n}S_j+\circ\left(\frac{1}{n}\right).$$ For the second term, we get
\begin{eqnarray*}
&&\mathbf{E}\left[\sum_{j=0}^{k-1}\left\Vert\int_{\frac{j}{n}}^{\frac{j+1}{n}}\int_{[0,1]}q(\theta_j]\mathbf{1}_{0\leq
x\leq
Re(Tr[\mathcal{J}(\theta_j)])})N(.,dx,ds)\right.\right.-\left(q(\overline{\rho}_j)+\circ\left(1\right)\right)\overline{\nu}_{j+1}(\overline{\rho}_j,.)]\Bigg{\Vert}\Bigg{]}\\
&=&
\mathbf{E}\left[\sum_{j=0}^{k-1}\left\Vert q(\theta_j)N(.,H_j(\theta_j))-\left(q(\overline{\rho}_j)+\circ\left(1\right)\right)\overline{\nu}_{j+1}
(\overline{\rho}_j,.)]\right\Vert\right]\\
&\leq&
\mathbf{E}\left[\sum_{j=0}^{k-1}\left\Vert q(\theta_j)\overline{\nu}_{j+1}(\overline{\rho}_j,.)-\left(q(\overline{\rho}_j)+\circ\left(1\right)\right)
\overline{\nu}_{j+1}(\overline{\rho}_j,.)]\right\Vert\right]\\
&&+\mathbf{E}\left[\sum_{j=0}^{k-1}\left\Vert q(\theta_j)\right\Vert\times\left\vert
N(.,H_j(\theta_j))-\overline{\nu}_{j+1}(\overline{\rho}_j,.)\right\vert\right]\\
&\leq&
\sum_{j=0}^{k-1}\mathbf{E}\left[(Q\Vert\theta_j-\overline{\rho_j}\Vert+\circ\left(1\right))\times\vert\overline{\nu}_{j+1}(\overline{\rho}_j,.)\vert\right]
\\&&+\sum_{j=0}^{k-1}\mathbf{E}[\Vert q(\theta_j)\Vert\times\vert
N(.,H_j(\theta_j))-\overline{\nu}_{j+1}(\overline{\rho}_j,.)\vert].
\end{eqnarray*}
Here we introduce the following discrete filtration \begin{equation}
\mathcal{F}_j=\sigma\left\{\overline{\nu}_{l}(\overline{\rho}_{l-1},.),N(.,H_l(\theta_l))/l\leq
j\right\}.
\end{equation}
It allows to compute the previous terms. It is clear that the random variables $\overline{\rho}_j$ and
$\theta_j$ are $\mathcal{F}_j$ measurable, we then have
 \begin{eqnarray*}
&&\mathbf{E}\left[\sum_{j=0}^{k-1}\left\Vert q(\theta_j)N(.,H_j(\theta_j))-\left(q(\overline{\rho}_j)+\circ\left(1\right)\right)
\overline{\nu}_{j+1}(\overline{\rho}_j,.)]\right\Vert\right]\\
&\leq&
\sum_{j=0}^{k-1}\mathbf{E}\left[(Q\Vert\theta_j-\overline{\rho_j}\Vert+\circ\left(1\right))\times\mathbf{E}
\left[\vert\overline{\nu}_{j+1}(\overline{\rho}_j,.)\vert/\mathcal{F}_j\right]\right]\\&&+\sum_{j=0}^{k-1}\mathbf{E}\left[\Vert
q( \theta_j)\Vert\times\mathbf{E}[\vert
N(.,H_j(\theta_j))-\overline{\nu}_{j+1}(\overline{\rho}_j,.)/\mathcal{F}_j\vert]\right].
\end{eqnarray*}
By conditioning with respect to $\mathcal{F}_j$, the random variable $\overline{\nu}_{j+1}(\overline{\rho}_j,.)$ is of Bernoulli type.
Hence we have
\begin{eqnarray*}\mathbf{E}[\vert\overline{\nu}_{j+1}(\overline{\rho}_j,.)\vert/\mathcal{F}_j]
&=&1-\exp(-\frac{1}{n}Tr(\mathcal{J}(\overline{\rho}_j))\\
&=&\frac{1}{n}Tr(\mathcal{J}(\overline{\rho}_j)+\circ\left(\frac{1}{n}\right).
\end{eqnarray*} For the second part, we have almost surely
\begin{eqnarray*}
&&\mathbf{E}\left[\vert
N(.,H_j(\theta_j))-\overline{\nu}_{j+1}(\overline{\rho}_j,.)/\mathcal{F}_j\vert\right]\\
&\leq& \mathbf{E}\left[\vert
N(.,H_j(\theta_j))-N(.,H_j(\overline{\rho}_j))/\mathcal{F}_j\right]+\mathbf{E}\left[N(.,H_j(\overline{\rho}_j))-\overline{\nu}_{j+1}(\overline{\rho}_j,.)/\mathcal{F}_j\right]\\
&\leq& \frac{1}{n}\left\vert
Tr[\mathcal{J}(\overline{\rho}_j))]-Tr[\mathcal{J}(\theta_j)]\right\vert+\mathbf{E}\left[N(.,H_j(\overline{\rho}_j))-\overline{\nu}_{j+1}(\overline{\rho}_j,.)/\mathcal{F}_j\right]\\
&\leq& \frac{1}{n}\left\vert
Tr[\mathcal{J}(\overline{\rho}_j))]-Tr[\mathcal{J}(\theta_j)]\right\vert+\left[\frac{1}{n}Tr[\mathcal{J}(\overline{\rho}_j)]-\left(1-\exp\left(-\frac{1}{n}Tr[\mathcal{J}(\overline{\rho}_j)]\right)\right)\right]\\
&\leq&
\frac{A}{n}\Vert\overline{\rho}_j-\theta_j\Vert+\frac{A'}{n^2}+\circ\left(\frac{1}{n^2}\right).
\end{eqnarray*}
The $\circ\left(\frac{1}{n^2}\right)$ are uniform in $j$ because $(\overline{\rho}_j)_j$ is uniformly bounded. For the second term,
the above inequalities and the fact that the Euler scheme is bounded implies that there exist two constants $K_1$ and $K_2$ such that
\begin{eqnarray*} &&\sum_{j=0}^{k-1}\mathbf{E}\left[\Vert q(
\theta_j)\Vert\times\mathbf{E}[\left\vert
N(.,H_j(\theta_j))-\overline{\nu}_{j+1}(\overline{\rho}_j,.)\right\vert/\mathcal{F}_j]\right]\\
&\leq&K_1\sum_{j=0}^{k-1}\frac{1}{n}S_j+\frac{K_2}{n}+\circ\left(\frac{1}{n}\right).
\end{eqnarray*}
For the first part, we have an equivalent inequality. Thus we can
conclude that there exist two constants $G_1$ and $G_2$ such
that
\begin{equation}
S_k\leq G_1\sum_{j=0}^{k-1}\frac{1}{n}S_j+\frac{G_2}{n}+\circ\left(\frac{1}{n}\right).
\end{equation}
The discrete Gronwall Lemma implies that there exists a constant
$F$ independent of $n$ such that for all $k\leq[nT]$
$$S_k\leq \frac{F}{n}.$$
The proposition is then proved.
\end{pf}\\

By using this two properties we can now express the final theorem.
\begin{thm}Let $T>0$ be a fixed time and let $(\Omega,\mathcal{F},\mathcal{F}_t,P)$ be the probability space of the poisson point process $N$.
Let $n$ be an integer and let $(\tilde{\rho}_{[nt]})_{0\leq t\leq T}$ be the discrete quantum trajectory defined for $k<[nT]$ by the equation
\begin{eqnarray*}
\tilde{\rho}_{k+1}&=&\mathcal{L}_{0}(\tilde{\rho}_k)+\mathcal{L}_{1}(\tilde{\rho}_k)\nonumber\\
&&+\left[-\frac{\mathcal{L}_{0}(\tilde{\rho}_k)}{Tr[\mathcal{L}_{0}(\tilde{\rho}_k)]}+
\frac{\mathcal{L}_{1}(\tilde{\rho}_k)}{Tr[\mathcal{L}_{1}(\tilde{\rho}_k)]}\right]
\left(\tilde{\nu}_{k+1}(\tilde{\rho}_k,.)-Tr[\mathcal{L}_{1}(\tilde{\rho}_k)]\right).
\end{eqnarray*}
Let $(\mu_t)_{0\leq t\leq T}$ be the quantum trajectory solution of the jump Belavkin equation on $[0,T]$ which satisfies
\begin{eqnarray*}
\mu_t&=&\mu_0+\int_0^tf(\mu_{s-})ds\\&&+\int_0^t\int_{[0,K]}\left[\frac{\mathcal{J}(\mu_{s-})}{Tr[\mathcal{J}(\mu_{s-})]}-\mu_{s-}\right]\mathbf{1}_{0\leq
x\leq Tr[\mathcal{J}(\mu_{s-})]}N(.,dx,ds).
\end{eqnarray*}

If $\mu_0=\tilde{ \rho}_0$, then the discrete quantum trajectory $(\rho_{[nt]})_{0\leq t\leq T}$ converges in distribution  to the continuous quantum trajectory $(\mu_t)_{0\leq t\leq T}$ in $\mathcal{D}\left([0,T]\right)$ for all $T$.
\end{thm}
\begin{pf}
Let $n$ be large enough. For $k\leq[nT]$, we define $\tilde{\mu}_k=\mu_{\frac{k}{n}}$. We define $$B_k=\mathbf{E}\left[\sup_{0\leq
i\leq k}\Vert \tilde{\rho}_i-\tilde{\mu}_i\Vert\right].$$ Thanks to Proposition $8$ and Theorem $5$ concerning the Euler scheme, there
exists a constant $R$ independent of $n$ such that for all $k\leq[nT]$
\begin{equation}\label{in}B_k\leq\frac{R}{n}.
\end{equation}

It is worth noticing that the process $(\tilde{\mu}_{[nt]})_{0\leq t\leq T}$ converges in distribution to $(\mu_t)_{0\leq t\leq T}$ for all $T$ in $\mathcal{D}\left([0,T]\right)$. Thanks to this fact and the inequality $(\ref{in})$, the convergence in distribution of $(\rho_{[nt]})_{0\leq t\leq T}$ to $(\mu_t)_{0\leq t\leq T}$ is proved.
\end{pf}

\end{document}